\documentclass[12pt,reqno]{amsart}
\usepackage{amssymb,amscd,amsthm}
\usepackage[T1]{fontenc}
\usepackage{amsmath,amssymb,graphicx,mathrsfs} 
\usepackage[colorlinks=true,allcolors = blue]{hyperref} 
\let\frak\mathfrak
\let\Bbb\mathbb

\def\>{\relax\ifmmode\mskip.666667\thinmuskip\relax\else\kern.111111em\fi}
\def\<{\relax\ifmmode\mskip-.333333\thinmuskip\relax\else\kern-.0555556em\fi}
\def\vsk#1>{\vskip#1\baselineskip}
\def\vv#1>{\vadjust{\vsk#1>}\ignorespaces}
\def\vvn#1>{\vadjust{\nobreak\vsk#1>\nobreak}\ignorespaces}

  \let\ssize\scriptstyle
\let\sssize\scriptscriptstyle

\let\Medskip\medskip
\def\medskip{\par\Medskip}
\let\Bigskip\bigskip
\def\bigskip{\par\Bigskip}

\let\Maketitle\maketitle
\def\maketitle{\Maketitle\thispagestyle{empty}\let\maketitle\empty}

\newtheorem{thm}{Theorem}[section]
\newtheorem{cor}[thm]{Corollary}
\newtheorem{lem}[thm]{Lemma}
\newtheorem{prop}[thm]{Proposition}

\numberwithin{equation}{section}

\theoremstyle{definition}
\newtheorem{defn}[thm]{Definition}
\newtheorem{rem}[thm]{Remark}
\newtheorem{prob}[thm]{Problem}

\let\mc\mathcal
\let\nc\newcommand

\let\al\alpha
\let\bt\beta
\let\dl\delta
\let\Dl\Delta
\let\eps\varepsilon
\let\gm\gamma
\let\Gm\Gamma

\let\la\lambda
\let\La\Lambda
\let\pho\phi
\let\phi\varphi
\let\si\sigma

\let\Om\Omega

\let\der\partial

\let\geq\geqslant

\let\leq\leqslant

\let\on\operatorname
\let\bi\bibitem
\let\bs\boldsymbol

\def\C{{\mathbb C}}
\def\Z{{\mathbb Z}}
\def\Q{{\mathbb Q}}

\def\Pb{{\mathbb P}}
\def\Qf{{\bf Q}}

\def\F{{\mc F}}

\def\+#1{^{\{#1\}}}

\def\diag{\on{diag}}

\def\Hom{\on{Hom}}

\def\tr{\on{tr}}

\def\sln{\mathfrak{sl}_N}

\def\beq{\begin{equation}}
\def\eeq{\end{equation}}
\def\be{\begin{equation*}}
\def\ee{\end{equation*}}

\nc{\bea}{\begin{eqnarray*}}
\nc{\eea}{\end{eqnarray*}}
\nc{\bean}{\begin{eqnarray}}
\nc{\eean}{\end{eqnarray}}
\nc{\bal}{\begin{align*}}
\nc{\eal}{\end{align*}}
\nc{\baln}{\begin{align}}
\nc{\ealn}{\end{align}}

\nc{\Il}{{\mc I_{\bs\la}}}
\nc{\bla}{{\bs\la}}
\nc{\Fla}{\F_\bla}
\nc{\tfl}{{T^*\Fla}}
\nc{\GL}{{GL_n(\C)}}
\nc{\GLC}{{GL_n(\C)\times\C^*}}

\let\sd s

\def\ddk_#1{\kk_{#1}\<\>\frac\der{\der\<\>\kk_{#1}}}

\def\bul{\mathbin{\raise.2ex\hbox{$\sssize\bullet$}}}
\def\intt{\mathchoice
{\mathop{\raise.2ex\rlap{$\,\,\ssize\backslash$}{\intop}}\nolimits}
{\mathop{\raise.3ex\rlap{$\,\sssize\backslash$}{\intop}}\nolimits}
{\mathop{\raise.1ex\rlap{$\sssize\>\backslash$}{\intop}}\nolimits}
{\mathop{\rlap{$\sssize\<\>\backslash$}{\intop}}\nolimits}}

\let\kk q
 
\let\cc c

\let\Ko K

\def\GZ/{Gelfand-Zetlin}
\def\KZ/{{\slshape KZ\/}}
\def\qKZ/{{\slshape qKZ\/}}
\def\XXX/{{\slshape XXX\/}}

\nc{\slnl}{{\sln (\lambda)}}
\nc{\PCN}{{   (\C[x])^N   }}
\nc{\di}{\on{Diag}}
\nc{\dio}{\on{Diag}_0}
\nc{\Mm}{{\mc M}}
\nc{\Nn}{{\mc N}}

\nc{\A}{{\mc C}}

\nc{\PCr}{{  P  (\C[x])^n   }}

\nc{\Pk}{{(\bs{P}^1)^k}}

\nc{\N}{{\Bbb N}}

\nc{\Ll}{{\mc L}}

\nc{\ord}{{\on{ord}\,}}

\nc{\Sing}{{\on{Sing}\,}}
\nc{\sing}{{\on{Sing}\,}}

\nc{\Hess}{{\on{Hess}}}

\nc{\R}{{\Bbb R}}

\let\on\operatorname
\nc{\Kk}{{\bs K}}
\nc{\Ap}{{A_\Phi(z)}}
\nc{\ap}{{A_\Phi(z)}}

\nc{\sv}{{\sing V}}
\nc{\cd}{{\C^n-\Delta}}
\nc{\UT}{{U^0}}   
\nc{\ep}{\epsilon}

\usepackage[all,cmtip]{xy}
\usepackage{empheq}
\usepackage{bm}

\newlanguage\fakelanguage
\newcommand\cyr{\fontencoding{OT2}\fontfamily{wncyr}\selectfont
   \language\fakelanguage}
\DeclareTextFontCommand{\textcyr}{\cyr}

\numberwithin{equation}{section}

\DeclareMathOperator{\HOM}{\mathscr{H}\text{\kern -3pt {\calligra\large om}}\,}

\makeatletter
\newsavebox{\@brx}
\newcommand{\llangle}[1][]{\savebox{\@brx}{\(\m@th{#1\langle}\)}%
  \mathopen{\copy\@brx\kern-0.5\wd\@brx\usebox{\@brx}}}
\newcommand{\rrangle}[1][]{\savebox{\@brx}{\(\m@th{#1\rangle}\)}%
  \mathclose{\copy\@brx\kern-0.5\wd\@brx\usebox{\@brx}}}
\makeatother

\newcommand{\bsh}{\begin{shaded}}
\newcommand{\esh}{\end{shaded}}

\newcommand{\dfm}{Dubrovin-Frobenius manifold}
\newcommand{\dfms}{\dfm s}
\newcommand{\egr}{{\rm gr}}

\begin{document}
\title[Degenerate RHB problems, Semisimplicity, and convergence]{Degenerate Riemann-Hilbert-Birkhoff problems, semisimplicity, and convergence of WDVV-potentials}
\author[Giordano Cotti ]{Giordano Cotti$\>^\circ$ }
\maketitle
\begin{center}

\textit{ $^\circ\>$Faculdade de Ci\^encias da Universidade de Lisboa\\ Grupo de F\'isica Matem\'atica \\
Campo Grande Edif\'icio C6, 1749-016 Lisboa, Portugal\/}

\end{center}
{\let\thefootnote\relax
\footnotetext{\vskip5pt 
\noindent
$^\circ\>$\textit{ E-mail}:  gcotti@fc.ul.pt, gcotti@sissa.it}}

\vskip6mm
\begin{abstract}
In the first part of this paper, we give a new analytical proof of a theorem of C.\,Sabbah on integrable deformations of meromorphic connections on $\Pb^1$. This theorem generalizes a previous result of B.\,Malgrange to the case of connections admitting irregular singularities of Poincar\'e rank 1 with coalescing eigenvalues.
In the second part of this paper, as an application, we prove that any semisimple formal Frobenius manifold (over $\C$), with unit and Euler field, is the completion of an analytic pointed germ of a \dfm. In other words, any formal power series, which provides a quasi-homogenous solution of WDVV equations and defines a semisimple Frobenius algebra at the origin, is actually convergent under no further tameness assumptions. 
\end{abstract}

\tableofcontents

\section{Introduction}

In this paper, we address the problem of convergence of formal solutions, in the ring of formal power series, of the Witten-Dijkgraaf-Verlinde-Verlinde (WDVV) associativity equations. This is the
overdetermined system of non-linear partial differential equations, in a single scalar function $F(t^1,\dots, t^n)$, given by
\begin{align*}
&\sum_{\mu,\nu}\frac{\der^3F}{\der t^\al\der t^\bt\der t^\mu}\eta^{\mu\nu}\frac{\der^3F}{\der t^\nu\der t^\gm\der t^\dl}=\sum_{\mu,\nu}\frac{\der^3F}{\der t^\dl\der t^\bt\der t^\mu}\eta^{\mu\nu}\frac{\der^3F}{\der t^\nu\der t^\gm\der t^\al},&\quad \al,\bt,\gm,\dl=1,\dots,n,\\
&\frac{\der^3F}{\der t^1\der t^\al\der t^\bt}=\eta_{\al\bt}={\rm const.,}\quad \eta=(\eta_{\al\bt})_{\al,\bt},\quad \eta^{-1}=(\eta^{\al\bt})_{\al,\bt}&\al,\bt=1,\dots, n.
\end{align*}
Introduced in the physics of topological field theories \cite{wit,DVV}, the geometry of solutions $F$ of WDVV equations, satisfying a further quasi-homogeneity condition
\[\sum_\al[(1-q_\al)t^\al+r^\al]\frac{\der F}{\der t^\al}=(3-d)F(\bm t)+\text{quadratic terms in }\bm t,
\]for suitable complex numbers $q_\al,r^\al,d\in\C$, was firstly axiomatized by B.\,Dubrovin, with the notion of Frobenius manifolds \cite{Dub92,Dub96,Dub98,Dub99}.
\vskip1.5mm
It was soon realized that these quasi-homogeneous solutions of WDVV equations arise in areas of mathematics which are very apart from each other (singularity theory, algebraic and symplectic geometry, integrable systems, mirror symmetry, to name just a few), often leading to new and non-trivial relations between them, see \cite{Dub96,man,Her02,Sab07}. 
\vskip1.5mm
Typically, the corresponding solutions $F(\bm t)$ of WDVV equations are given as generating functions of numerical sequences of geometrical interest (e.g. Gromov-Witten invariants). Consequently, they can be handled just as formal power series in $k[\![\bm t]\!]$, where $k$ is a commutative $\Q$-algebra, defining a \emph{formal Frobenius manifold} structure on the formal spectrum ${\rm Spf }\,k[\![\bm t]\!]$, see \cite[III.\S1]{man}. This defines a formal family of Frobenius algebras with structure constants given by $c_{\al\bt}^\gm(\bm t):=\eta^{\la\gm}\der_{\al\bt\la}^3F(\bm t)$.
\vskip1.5mm 
The relevance of these formal structures is further highlighted by their deep relations with the cohomology of the Deligne-Mumford moduli stacks $\overline{\mc M}_{g,\frak n}$ of $\frak n$-pointed stable curves of genus $g$, \cite{KM94,man}. Remarkably enough, any formal Frobenius manifold is equivalent to a {\it tree level\footnote{A richer notion of {\it complete} CohFT on a given $(H,\eta)$ is also available, in which the datum is enriched to a family $(\Om_{g,\frak n})_{g,\frak n}$ of $k$-linear tensors $\Om_{g,\frak n}\in (H^*)^{\otimes \frak n}\otimes_k H^\bullet(\overline{\mc M}_{g,\frak n}; k)$, satisfying further compatibility properties, for any pair $(g,\frak n)$ of non-negative integers in the stable range $2g-2+\frak n>0$. The prototypical example of a complete CohFT is provided by the Gromov-Witten theory of a smooth projective variety $X$. The corresponding formal Frobenius manifold attached to its genus zero sector is called \emph{quantum cohomology of $X$}. See \cite{KM94,man} and Section \ref{sec6} of this paper.} Cohomological Field Theory} (CohFT), i.e. the datum of a family of $\frak S_\frak n$-covariant\footnote{Here $\frak S_n$ denotes the symmetric group on a finite set with $n$ elements.} tensors $\Om_{0,\frak n}\in (H^*)^{\otimes \frak n}\otimes_k H^\bullet(\overline{\mc M}_{0,\frak n}; k)$, on a given free metric $k$-module $(H,\eta)$ of finite rank, satisfying some compatibility conditions with respect to the natural forgetful morphisms $\overline{\mc M}_{0,\frak n}\to \overline{\mc M}_{0,\frak n-1}$, and gluing morphisms $\overline{\mc M}_{0,\frak n_1}\times \overline{\mc M}_{0,\frak n_2}\to \overline{\mc M}_{0,\frak n_1+\frak n_2}$.  
The corresponding WDVV-potential $F(\bm t)$ is a generating power series for integrals of the form $\int_{\overline{\mc M}_{0,\frak n}}\Om_{0,\frak n}(\bigotimes_{j=1}^\frak n \Dl_{\al_j})$ for a $k$-basis $(\Dl_j)_{j}$ of $H$. See \cite{KM94,man,Pan18} for more details.

\vskip1.5mm
One of the main points of the current paper is to find sufficient conditions ensuring the convergence of quasi-homogeneous solutions $F\in k[\![\bm t]\!]$ of WDVV equations, in the case $k=\C$. 
The convergence condition allows to jump from the formal category to the complex analytic category: formal Frobenius manifolds can be promoted to \emph{analytic} Frobenius manifolds, the class of geometrical objects originally conceived by Dubrovin, and for this reason also called {\it \dfm s}.
\vskip1.5mm
The new main result of this paper, Theorem \ref{thconv}, claims that any formal {\it semisimple} Frobenius manifold over $k=\C$ is actually the completion of a pointed germ of an analytic \dfm. Alternatively stated, given a quasi-homogeneous formal solution $F\in\C[\![\bm t]\!]$ of WDVV equations whose corresponding Frobenius $\C$-algebra at the origin $\bm t=0$ is semisimple, its domain of convergence is non-empty, and it thus carries a \dfm\ structure. This statement is a refinement of a seemingly known result, referred to as a ``general fact'' in \cite[III.\S7.1, pag.135]{man}, and stated under stronger unnecessary \emph{tameness} assumptions\footnote{I do not know any reference in literature where a complete proof is given. I thank Yu.I.\,Manin for a friendly e-mail correspondence on this point. The current paper both recovers a proof of this known fact, and it also removes the tameness assumption.} (see the next paragraph).
\vskip1.5mm
 At the core of our proof there is the local identification of semisimple points $\bm t$ of a \dfm\ with the parameters of isomonodromic deformations of ordinary differential equations with rational coefficients, of the form
\beq\label{ODE}
\frac{d}{dz}Y(z,\bm t)=\left(\mc U(\bm t)+\frac{1}{z}\mu(\bm t)\right)Y(z,\bm t),
\eeq
where $\mc U,\mu$ are (matrices representing) suitably defined tensors on the \dfm. This identification -- one of the main points of the theory of Dubrovin-- was originally established in \cite{Dub96,Dub98,Dub99} at {\it tame semisimple points}, i.e. points $\bm t$ at which the leading term $\mc U(\bm t)$ of the coefficient of \eqref{ODE} has simple spectrum. Subsequently, in \cite{CG17,CG18,CDG1,CDG20} the isomonodromic approach to the Frobenius geometry was extended to {\it all} semisimple points, including points $\bm t$ at which some of the eigenvalues of $\mc U(\bm t)$ coalesce.
\vskip1.5mm
The proof of Theorem \ref{thconv} consists of two parts. Firstly, given a formal Frobenius manifold $F\in\C[\![\bm t]\!]$, it is constructed an \emph{analytic family} \eqref{ODE} of ODEs specializing to the given one\footnote{Given a formal Frobenius manifold, the system \eqref{ODE} has coefficients in $M_n(\C[\![\bm t]\!])$. Hence, for $t=0$, we have a well defined differential system with coefficients in $M_n(\C)$.} for $\bm t=0$, and defining a \dfm. Secondly, it is proved that the underlying analytic\footnote{Here $\C\{\bm t\}$ denotes the algebra of convergent power series in $\bm t$.} WDVV-potential $F^{\rm an}\in\C\{\bm t\}$ coincides with the original formal one, i.e. $F=F^{\rm an}$. 
It is thus clear that the first step of the proof of Theorem \ref{thconv} relies on the existence of solutions of \emph{families} of Riemann-Hilbert-Birkhoff boundary value problems. In the case $\bm t=0$ is a tame semisimple point of the given formal Frobenius manifold, a well-known result of B.\,Malgrange \cite{Mal83a,Mal83b,Mal86}, on the existence of universal integrable deformations of meromorphic connections on $\Pb^1$ with irregular singularities, can be applied. This leads to the already known result mentioned\footnote{We warn the reader that in the exposition of \cite{man}, the isomonodromic system \eqref{ODE} is replaced by a Fuchsian one obtained by applying a (formal) Laplace transform, see \cite[Ch. II.\S1-3]{man}.} in \cite[III.\S7, pag.135]{man}. 
\vskip1.5mm
In \cite{Sab18}, C.\,Sabbah obtained an extension of the theorem of Malgrange, in order to include the case of meromorphic connections on $\Pb^1$ which admit irregular singularities with coalescing eigenvalues. In the geometrical case attached to Frobenius manifolds, the assumptions of \cite[Th. 4.9]{Sab18} are satisfied. Sabbah Theorem can thus be applied in the first step of the proof of Theorem \ref{thconv}, in the case $\bm t=0$ is a coalescing semisimple point for the given formal Frobenius manifold. Remarkably enough, the assumptions of \cite[Th. 4.9]{Sab18} exactly coincide with the sharp conditions, found in \cite{CG17,CDG1}, under which the resulting analytic family \eqref{ODE} of ODEs has a well-behaved deformation theory of both formal and genuine solutions.
\vskip1.5mm
The original proof of \cite[Th. 4.9]{Sab18} is actually only one of the outcomes of a more general study, invoking a mix of techniques, including properties of  \emph{good} and \emph{very-good formal decompositions} of flat meromorphic bundles \cite{Sab93,Sab00}, and recent results on meromorphic connections in dimension $\geq 2$ due to K.\,Kedlaya (in the complex analytic case) and T.\,Mochizuki (in the algebraic case), see \cite{Ked10,Ked11,Moc09,Moc11a,Moc11b,Moc14}. In Section \ref{secsabth}, we give an alternative proof of \cite[Th. 4.9]{Sab18}, with a more analytical perspective, closer to the one of \cite{CG17,CDG1}. Our proof is uniquely based on properties of Fredholm-operator-valued holomorphic functions. In particular, a result due to B.\,Gramsch \cite{Gra70} -- an {\it analytical Fredholm alternative} with respect to several parameters--  will be invoked to prove that the solvability of a family of Riemann-Hilbert-Birkhoff boundary value problems is an open property, in the same spirit of \cite{Zho89}. This is a well-known strategy for proving the Painlev\'e property of solutions of the isomonodromy deformation equations, see e.g. \cite{FZ92,Pain}. 
\vskip1.5mm
Many of the results of this paper can be extended to the case of flat $F$-manifolds \cite{HM99,Man05}. These are slightly weaker structures than the Frobenius one, but whose geometry encompasses even more areas of modern mathematics, such as special solutions of the oriented associativity equations \cite{LM04}, quantum $K$-theory \cite{Lee04}, all Painlev\'e transcendents \cite{AL15}, open WDVV equations \cite{BB19}, $F$-cohomological field theories \cite{ABLR20}, and even information geometry \cite{CM20}. We plan to give more analytical details in a future publication.
\vskip1.5mm
\noindent{\bf Structure of the paper. }In Section \ref{sec2} we review necessary background material on the Riemann-Hilbert-Birkhoff problems with a geometrical perspective. The main results of B.\,Malgrange on the existence of the universal integrable deformation of meromorphic connection, as well as their generalization to degenerate cases due to C.\,Sabbah, are presented and summarized. 
\vskip1mm
Section \ref{secsabth} is devoted to an analytical proof of Malgrange-Sabbah Theorem. After introducing the notion of admissible data, we formulate a Riemann-Hilbert-Birkhoff boundary value problem $\mc P[\bm u,\tau,\frak M]$, depending on parameters $\bm u\in\C^n$. We factorize its solutions via two auxiliary RHB problems, and we analyze its solvability with respect to $\bm u$. 
\vskip1mm
In Section \ref{sec4} basic notions in the theory of both formal and analytic Frobenius manifolds are given. We explain how to pass from the analytic to the formal category, and vice-versa under convergence assumption of the WDVV potential.
\vskip1mm
In Section \ref{sec5}, we review necessary results on the extended deformed connection on both formal and analytic Frobenius manifolds, properties of solutions of the Darboux-Egoroff system of partial differential equations, and the reconstruction procedure of the Frobenius potential. Consequently, we prove the main new result, Theorem \ref{thconv}. 
\vskip1mm
In the last Section \ref{sec6}, reformulations and applications of Theorem \ref{thconv} to cohomological field theories and quantum cohomology are given.
\vskip1.5mm
\noindent{\bf Acknowledgements. }The author is thankful to P.\,Belmans, M.\,Bertola, G.\,Bogo, G.\,Broc\-chi, G.M. Dall'Ara, I.\,Gayur, T.\,Grava, D.\,Guzzetti, C.\,Hertling, A.R.\,Its, P.\,Lorenzoni, Yu.I. Man\-in, D.\,Masoero, M.\,Maz\-zocco, A.T.\,Ricolfi, P.\,Rossi, V.\,Roubtsov, A.\,Var\-chenko, C. Sabbah, M. Smirnov, D.\,Yang for several valuable discussions. The author is thankful to the Hausdorff Research Institute for Mathematics (HIM) in Bonn, Germany, where this project was completed, for providing excellent working conditions during the JTP ``New Trends in Representation Theory''. The author is also grateful to the anonymous Referees, for helpful comments and suggestions improving the exposition. This research was supported by the EPSRC Research Grant EP/P021913/2, by HIM (Bonn, Germany), and by the FCT Project PTDC/MAT-PUR/ 30234/2017 ``Irregular connections on algebraic curves and Quantum Field Theory''.

\section{Degenerations of Riemann-Hilbert-Birkhoff inverse problems }\label{sec2}
\subsection{Riemann-Hilbert-Birkhoff inverse problems}Let $D$ be a disc centered at $z=\infty$ in $\Pb^1$. Given a (trivial) vector bundle on $D$ equipped with a meromorphic connection with a pole at $z=\infty$, the Riemann-Hilbert-Birkhoff (RHB) problem is the following:

\begin{prob}\label{RHB1}
Does there exist a {\it trivial} vector bundle $E^o$ on $\Pb^1$ equipped with a meromorphic connection $\nabla^o$, restricting to the given data on $D$, and with a further {\it logarithmic} pole only at $z=0$? 
\end{prob}

Assume that the pole at $z=\infty$ is of Poincar\'e rank 1: in a basis of sections on $D$, the meromorphic connection has matrix of connection 1-forms $\Om=-A(z)dz$, where the $n\times n$ matrix $A(z)$ equals
\[A(z)=\sum_{k=0}^\infty A_kz^{-k},\quad A_0\neq0.
\]Denote by $\mc O\{\frac{1}{z}\}$ the ring of convergent power series in $\frac{1}{z}$. The RHB problem \ref{RHB1} is then equivalent to find a so-called {\it Birkhoff normal form}: does it exist a matrix $G\in GL_n(\mc O\{\frac{1}{z}\})$ such that $B(z)=G^{-1}AG-G^{-1}\frac{d}{dz}G$ is of the form
\[B(z)=B_0+\frac{B_1}{z},\quad B_0,B_1\in M_n(\C)?
\]
\subsection{Universal integrable deformations of Birkhoff normal forms: Malgrange Theorems} 
In this paper we consider {\it families} of RHB problems parametrized by a parameter space $X$, see \cite{Mal83a,Mal83b,Mal86}\cite[Ch.VI]{Sab07}. We assume $X$ is a complex manifold.

\begin{defn}
Let $(E^o,\nabla^o)$ be a trivial vector bundle on a disc $D$ equipped with a meromorphic connection with a pole of order 2 at $z=\infty$. An {\it integrable deformation} of $(E^o,\nabla^o)$ parametrized by $X$ is the datum $(E,\nabla)$ of 
\begin{itemize}
\item a trivial vector bundle $E$ on $D\times X$,
\item a {\it flat} connection $\nabla$ on $E$ with a pole of Poincar\'e rank 1 along $\{\infty\}\times X$,
\end{itemize}
such that $(E,\nabla)$ restricts to $(E^o,\nabla^o)$ at a point $x_o\in X$. The integrable deformation is called {\it versal} if any other deformation with base space $X'$ is induced by the previous one via pull-back by a holomorphic map $\phi\colon (X',x_o')\to (X,x_o)$. It is {\it universal} if the germ at $x_o'$ of the base-change $\phi$ is uniquely determined.
\end{defn}

Let $(E^o,\nabla^o)$ be a solution of a RHB problem \ref{RHB1}, i.e. a trivial vector bundle on $\Pb^1$ with meromorphic connection with matrix (in a suitable basis of sections) of the form
\beq\label{con1}\Om=-\left(\La_o+\frac{B_o}{z}\right)dz.
\eeq
\vskip2mm
Recall that a matrix $A\in M_n(\C)$ is said to be {\it regular} if any (and hence all) of the following equivalent conditions is satisfied:
\begin{enumerate}
\item the characteristic polynomial of $A$ equals its minimal polynomial,
\item the commutator of $A$ in $M_n(\C)$ is of minimal dimension (i.e. it equals $n$),
\item the commutator of $A$ in $M_n(\C)$ is $\C[A]$.
\end{enumerate}

\begin{thm}[{\cite{Mal83a,Mal86}}]Assume that the matrix $\La_o$ is regular. The connection $\nabla^o$  with matrix \eqref{con1} has a germ of universal deformation.
\end{thm}

This result can be refined to a global statement, under the further {\it semisimplicity} assumption on $\La_o$. Let us then assume that $\La_o={\rm diag}(u^1_o,\dots, u^n_o)$ with $u_o^i\neq u_o^j$ for $i\neq j$.
\vskip2mm
Let $\Dl$ be the union of {big diagonal} hyperplanes in $\C^n$, defined by the equations
\[\Dl:=\bigcup_{i<j}\{\bm u\in\C^n\colon u^i=u^j\},
\]let $X_n$ be the complement $\C^n\setminus \Dl$, with base point $\bm u_o:=(u_o^1,\dots, u_o^n)$. Denote by $\pi\colon (\widetilde X_n,\tilde{ \bm u}_o)\to (X_n,\bm u_o)$ the universal cover of $X_n$, equipped with fixed base points $\tilde {\bm u}_o$ and $\bm u_o$, respectively. The space $X_n$ is identified with the space of diagonal regular $n\times n$ matrices.

\begin{thm}[\cite{JMU81,Mal83b}]\label{tmal}
There exists on $\Pb^1\times \widetilde X_n$ a vector bundle $E$, equipped with a meromorphic connection $\nabla$, such that
\begin{enumerate}
\item  $\nabla$ is flat, with a pole of Poincar\'e rank 1 along $\{\infty\}\times \widetilde X_n$, and a logarithmic pole along $\{0\}\times \widetilde X_n$, 
\item $(E,\nabla)$ restricts to $(E^o,\nabla^o)$ at $\tilde{\bm u}_o$,
\item for any $\tilde{\bm u}\in\widetilde X_n$, the eigenvalues of the residue of $\nabla$ at the point $(\infty,\tilde{\bm u})$ equal (up to permutation) the $n$-tuple $\pi(\tilde{\bm u})$.
\end{enumerate}
Let $\Theta\subseteq \widetilde X_n$ be the hypersurface of points $\tilde{\bm u}\in\widetilde X_n$ such that $E|_{\Pb^1\times\{\tilde{\bm u}\}}$ is not trivial. The coefficients of $\nabla$ have poles along $\Theta$. Moreover, for any $\tilde{\bm u}\in\widetilde X_n\setminus\Theta$, the bundle with meromorphic connection $(E,\nabla)$ induces a universal deformation of its restriction $(E,\nabla)|_{\Pb^1\times\{\tilde{\bm u}\}}$.
\end{thm}

It is possible to explicitly describe the matrix of the connection 1-forms of the universal deformation of Theorem \ref{tmal}. 
\vskip2mm
For $\bm u\in\C^n$, denote $\La(\bm u):=\diag(u^1,\dots,u^n)$, so that $\La(\bm u_o)=\La_o$. Given a matrix $A$ denote by $A'$ its diagonal part, and by $A''$ its off-diagonal part.
\vskip2mm
For $\bm u_o\notin\Dl$, there exists an off-diagonal matrix $F''(\bm u)$, holomorphic near $\bm u_o$, such that the flat connection $\nabla$ of Theorem \ref{tmal} has matrix of connection 1-forms
\beq\label{uic}
-d\left(z\La(\bm u)\right)-\left([\La(u),F''(\bm u)]+B'_o\right)\frac{dz}{z}-[d\La(\bm u),F''(\bm u)],
\eeq
e.g. see \cite[VI.\S 3.f, eq. (3.12)]{Sab07}. 
Notice that the $dz$-component of \eqref{uic} restricts to \eqref{con1} at $\bm u=\bm u_o$. Moreover, there exists a $z^{-1}$-formal base change which transforms \eqref{uic} into
\beq\label{con2}
-d\left(z\La(\bm u)\right)-B'_o\frac{dz}{z}.
\eeq

\subsection{Integrable deformations of degenerate Birkhoff normal forms: Sabbah Theorem}  In the notations of the previous section, assume $\bm u_o\in\Dl$. Define the partition $\{1,\dots,n\}=\coprod_{r\in R}I_r$ such that for any $r\in R$ we have
\[\{i,j\}\subseteq I_r\quad\text{if and only if}\quad u_o^i=u_o^j.
\]  In \cite{Sab18}, C.\,Sabbah addressed the following problem. 
\vskip2mm
{\bf Question: }Is it possible to find an integrable deformation of the form \eqref{uic} of the Birkhoff normal form \eqref{con1} with $z^{-1}$-formal normal form \eqref{con2}? 
\vskip2mm
Remarkably, in \cite[Section 4]{Sab18} it is shown that the answer is positive, under (sharp) sufficient conditions on the coefficient $B_o$ of the normal form \eqref{con1}.

\begin{thm}[{\cite[Th. 4.9]{Sab18}}]\label{tsab}
Let $\bm u_o\in\Dl$, and $\mc V$ a neighborhood of $\bm u_o$ in $\C^n$. Assume that
\begin{enumerate}
\item $B''_o\in {\rm Im\ ad}(\La(\bm u^o))$,
\item $B'_o$ is \emph{partially non-resonant}, i.e. 
\[\forall\ r\in R,\quad \forall\ i,j\in I_r,\quad (B'_o)_{ii}-(B'_o)_{jj}\notin \Z\setminus\{0\}.\]
\end{enumerate}
If $\mc V$ is sufficiently small, there exists a holomorphic hypersurface $\Theta$ in $\mc V\setminus\{\bm u_o\}$ and a holomorphic off-diagonal matrix $F''(\bm u)$ on $\mc V\setminus\Theta$, such that the meromorphic connection, on the trivial vector bundle on $\mathbb P^1\times (\mc V\setminus\Theta)$, with matrix \eqref{uic} is integrable, restricts to \eqref{con1} at $\bm u_o$, and it is formally equivalent at $z=\infty$ to the matrix connection \eqref{con2}.
\end{thm}

\section{An analytical proof of Sabbah Theorem}\label{secsabth}
In this section we provide an analytical proof of Sabbah Theorem \ref{tsab}, based on properties of holomorphic Fredholm-operator-valued functions. We recast Sabbah Theorem in terms of suitable Riemann-Hilbert-Birkhoff boundary value problems $\mc P[\bm u,\tau,\frak M]$ depending on parameters $(\bm u,\tau,\frak M)$, the {\it admissible data}. Sabbah Theorem claims that the solvability of $\mc P[\bm u,\tau,\frak M]$, for fixed $(\tau,\frak M)$, is an open property with respect to $\bm u$. We prove this statement by factorizing solutions of $\mc P[\bm u,\tau,\frak M]$ in terms of two auxiliary RHB boundary value problems, $\mc P_1[\bm u,\tau,\frak M]$ and $\mc P_2[\bm u,\tau,\frak M]$. Their solvability  is studied via results in the $L^p$-theory of RH boundary value problems. General references for this section are \cite{AB,Bot20,CG81,CG18,CDG1,DZ02a,DZ02b,Pain,Its03,Its11,MP80,TO16,Vek67,Zho89}.
\subsection{Admissible data and Riemann-Hilbert-Birkhoff boundary value problem} Denote by ${\rm Arg}(z)\in ]-\pi,\pi]$ the principal branch of the argument of the complex number $z$.  Let $\bm u\in\C^n$, and set
\[
\mathscr S(\bm u):=\left\{{\rm Arg}\left[-\sqrt{-1}(\overline{u^i}-\overline{u^j})\right]+2\pi k\colon k\in\Z,\ i,j\text{ are s.t. }u^i\neq u^j\right\}.
\]
Any element $\tau\in\R\setminus \mathscr S(\bm u)$ will be said to be \emph{admissible at $\bm u$}.

\begin{defn}\label{datum}
Let $\bm u\in\C^n$ and $\tau$ be admissible at $\bm u$. A $(\bm u,\tau)$-{\it admissible datum} is a $6$-tuple $\frak M:=(B,D,L,S_1,S_2,C)$ of matrices in $M_n(\C)$ such that:
\begin{enumerate} 
\item the matrix $B$ is diagonal, i.e. $B=B'$,
\item $D$ is a diagonal matrix of integers,
\item we have
\beq
\label{vad1}
\tr B=\tr D+\tr L.
\eeq
\item the matrices $S_1,S_2,C$ are invertible, with $\det S_1=\det S_2=1$,
\item $(S_1)_{ii}=(S_2)_{ii}=1$, 
\item if $i\neq j$, then $(S_1^{-1})_{ij}=0$ if ${\rm Re}\left(e^{\sqrt{-1}(\tau-\pi)}(u^i-u^j)\right)>0$, 
\item if $i\neq j$, then $(S_2)_{ij}=0$ if ${\rm Re}\left(e^{\sqrt{-1}\tau}(u^i-u^j)\right)>0$,
\item we have
\beq\label{c2} S_1^{-1}e^{2\pi\sqrt{-1}B}S_2^{-1}=C^{-1}e^{2\pi\sqrt{-1}L}C.
\eeq
\end{enumerate}
If $\bm u\in\Dl$, define the partition $\{1,\dots,n\}=\coprod_{r\in R}I_r$ such that for any $r\in R$ we have
$\{i,j\}\subseteq I_r\quad\text{if and only if}\quad u^i=u^j$. We then require the further vanishing condition 
\begin{enumerate}
\item[(9)] $(S_1^{-1})_{ij}=(S_2)_{ij}=0$ if $i,j\in I_r$, with $i\neq j$, for some $r\in R$.
\end{enumerate}
\end{defn}

\begin{lem}
Let $\bm u_o\in \C^n$ and $\tau$ admissible at $\bm u_o$. If $\frak M$ is $(\bm u_o,\tau)$-admissible, then there exists a sufficiently small neighborhood $\mc V$ of $\bm u_o$ such that
\begin{enumerate}
\item $\tau$ is admissible at $\bm u$, for all $\bm u\in\mc V$,
\item $\frak M$ is $(\bm u,\tau)$-admissible for all $\bm u\in\mc V$. \qed
\end{enumerate}
\end{lem}

\begin{figure}
\centering
\def\svgscale{.6}
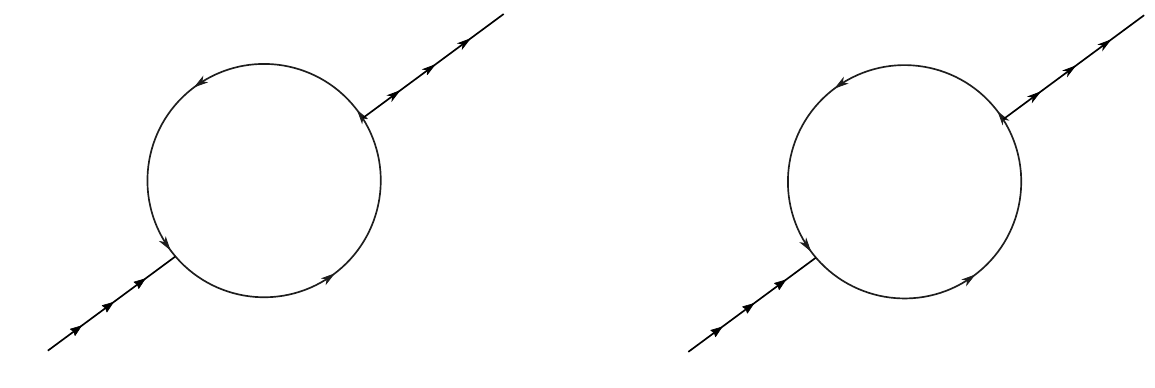
\caption{Contour $\Gm$, paths $\Gm_{\pm\infty},\Gm_1,\Gm_2$, domains $\Pi_0,\Pi_L,\Pi_R$, and $\pm$ sides of $\Gm$.}
\label{disegno1}
\end{figure}

Let $\bm u\in\C^n$ and $\tau$ admissible at $\bm u$. Consider the complex $z$-plane with a branch cut from $0$ to $\infty$:
\[\tau-\pi<\arg z<\tau+\pi.
\]
Let $r>0$ and denote by $\Gm=\Gm(\tau,r)$ the union of the following oriented paths, see Figure \ref{disegno1}:
\begin{enumerate}
\item the half-line $\Gm_{-\infty}$ defined by $\arg z=\tau\pm\pi$, $|z|>r$, originating from $\infty$;
\item the half-line $\Gm_{+\infty}$ defined by $\arg z=\tau$, $|z|>r$, ending to $\infty$;
\item the half-circle $\Gm_1$ defined by $\tau-\pi<\arg z<\tau$, $|z|=r$, counterclockwise oriented;
\item the half-circle $\Gm_2$ defined by $\tau<\arg z<\tau+\pi$, $|z|=r$, counterclockwise oriented.
\end{enumerate}
The orientations uniquely define the + and - side for each path $\Gm_{\pm\infty},\Gm_1,\Gm_2$. For $z\in\Gm_{-\infty}$ we use the symbol $z_\pm$ if $\arg z=\tau\pm\pi$. Set $\Pi_0,\Pi_L,\Pi_R$ to be the components of complement $\C\setminus\Gm$, and $T_1, T_2$ to be the two nodes of $\Gm$, as in Figure \ref{disegno1}.
\vskip2mm
Let $\frak M:=(B,D,L,S_1,S_2,C)$ be a $(\bm u,\tau)$-admissible datum. Define two functions
\[Q(-;\bm u), H(-;\bm u)\colon \Gm\to GL(n,\C),
\] by
\[Q(z;\bm u):=\La(\bm u)z+B\log z,\quad \La(\bm u):={\rm diag}(u^1,\dots, u^n),
\]
\begin{equation}\label{defH}
H(z;\bm u){:=}\left\{
\begin{alignedat}{2}
&e^{Q(z_-;\bm u)}S_1^{-1} e^{-Q(z_-;\bm u)},&\text{ along }\Gm_{-\infty},\\
       &e^{Q(z;\bm u)}S_2 e^{-Q(z;\bm u)},&\text{ along }\Gm_{+\infty},\\
       &e^{Q(z;\bm u)} C^{-1}z^{-L}z^{-D},&\text{ along }\Gm_{1},\\
       &e^{Q(z;\bm u)} S_2^{-1}C^{-1}z^{-L}z^{-D},&\text{ along }\Gm_{2}.
\end{alignedat}
\right.
\end{equation}

We denote by $H_{\pm\infty},H_1,H_2$ the restrictions of $H$ at $\Gm_{\pm\infty},\Gm_1,\Gm_2$.   These functions can be analytically continued, with respect to $z$, to the universal cover $\widetilde{\C^*}$. 

\begin{prop} 
The following identities hold true identically in $z\in\widetilde{\C^*}$:
\begin{align}
\label{sc1}H_{-\infty}(z;\bm u)H_2(z e^{2\pi\sqrt{-1}};\bm u)H_1(z;\bm u)^{-1}=I,\\
\label{sc2}H_1(z;\bm u)H_2(z;\bm u)^{-1}H_{+\infty}(z;\bm u)^{-1}=I.
\end{align}
\end{prop}
\proof
A simple computation shows that \eqref{sc1} follows from \eqref{c2}. Equation \eqref{sc2} is easily checked.
\endproof
    
\begin{prob}[{Problem $\mc P[\bm u,\tau,\frak M]$}] Find an analytic function $G\colon \C\setminus\Gm\to M_n(\C)$ such that 
\begin{enumerate}
\item $G|_{\Pi_\nu}$ extends continuously to $\overline{\Pi_\nu}$ for $\nu=0,L,R$;
\item the non-tangential limits $G_\pm\colon \Gm\setminus\{T_1,T_2\}\to M_n(\C)$ of $G$, from the $-$ and $+$ sides of $\Gm$, are related by
\[G_+(z)=G_-(z)H(z;\bm u);
\]
\item $G(z)$ tends to the identity matrix $I$ as $z\to\infty$.
\end{enumerate}
\end{prob}

\begin{rem}\label{26.05.21-rem1}
Let $G$ be a solution of $\mc P[\bm u,\tau,\frak M]$. Let $\widetilde{\C^*}$ be the universal cover of $\C^*$. 
Define the functions $Y_i(-;\bm u)\colon \widetilde{\C^*}\to \C$, with $i=0,1,2,3$, by
\bea
&Y_{0}(z;\bm u):=G(z;\bm u)z^{\mc D}z^L,\quad &z\in\Pi_{0},\\
& Y_1(z;\bm u):= G(e^{2\pi\sqrt{-1}}z;\bm u)z^{B'_o}e^{\La(\bm u)z},\quad &z\in e^{-2\pi\sqrt{-1}}\Pi_L,\\
&Y_{2}(z;\bm u):=G(z;\bm u)z^{B'_o}e^{\La(\bm u)z},\quad &z\in\Pi_{R},\\
&Y_{3}(z;\bm u):=G(z;\bm u)z^{B'_o}e^{\La(\bm u)z},\quad &z\in\Pi_{L}.
\eea
We have
\beq\label{26.05.21-1}
Y_2(z;\bm u)=Y_1(z;\bm u)\,S_1,\quad Y_3(z;\bm u)=Y_2(z;\bm u)\,S_2,\quad 
Y_2(z;\bm u)=Y_0(z;\bm u)\,C.
\eeq
It follows that
\[
\frac{\der Y_0}{\der z}Y_0^{-1}=\frac{\der Y_1}{\der z}Y_1^{-1}=\frac{\der Y_2}{\der z}Y_2^{-1}=\frac{\der Y_3}{\der z}Y_3^{-1},
\]
and the resulting function $\mc A(-;\bm u):=\der_z Y_i(-;\bm u)\cdot Y_i(-;\bm u)^{-1}$ is analytic with respect to $z\in\C^*$. Hence, for any fixed $\bm u$, the matrices $Y_i(-;\bm u)$, with $i=0,1,2,3$, are fundamental systems of solutions of the differential system
\beq\label{26.05.21-2}
\frac{dY}{dz}=\mc A(z;\bm u)Y.
\eeq
Such a system will be studied in Section \ref{ssabt}. The admissible data $\frak M$ determine the monodromy of the solutions $Y_i$. The matrices $S_1,S_2$ are the {\it Stokes matrices} of the system, the matrix $B$ is the {\it formal monodromy}, and the matrix $C$ is the {\it central connection matrix}. The independence of $\frak M$ with respect to $\bm u$ implies that the system \eqref{26.05.21-2} is {\it isomonodromic}.
Equations \eqref{26.05.21-1} are the reason for the precise shape of the RHB problem $\mc P[\bm u,\tau,\frak M]$.
\end{rem}

\subsection{Factorization of solutions}We factorize solutions of the problem $\mc P[\bm u,\tau,\frak M]$ via two auxiliary RHB boundary value problems, $P_1[\bm u,\tau,\frak M]$ and $\mc P_2[\bm u,\tau,\frak M]$. First, we describe the contours for both problems.
\vskip3mm
Let $P_1\in\Gm_{-\infty}$, and $P_2\in\Gm_{+\infty}$ such that $|P_1|,|P_2|>r$. Set 
\begin{itemize}
\item $\ell_1\subseteq \Gm_{-\infty}$ to be the closed half-line from $\infty$ to $P_1$, 
\item $\ell_2\subseteq \Gm_{+\infty}$ to be the closed half-line from $P_2$ to $\infty$.
\end{itemize}
Define 
\begin{itemize}
\item $\Gm'$ to be the union $\ell_1\cup \ell_2$,
\item $\Gm''$ to be a circle of radius $R>\max\{|P_1|,|P_2|\}$.
\end{itemize}

\begin{figure}
\centering
\def\svgscale{.6}
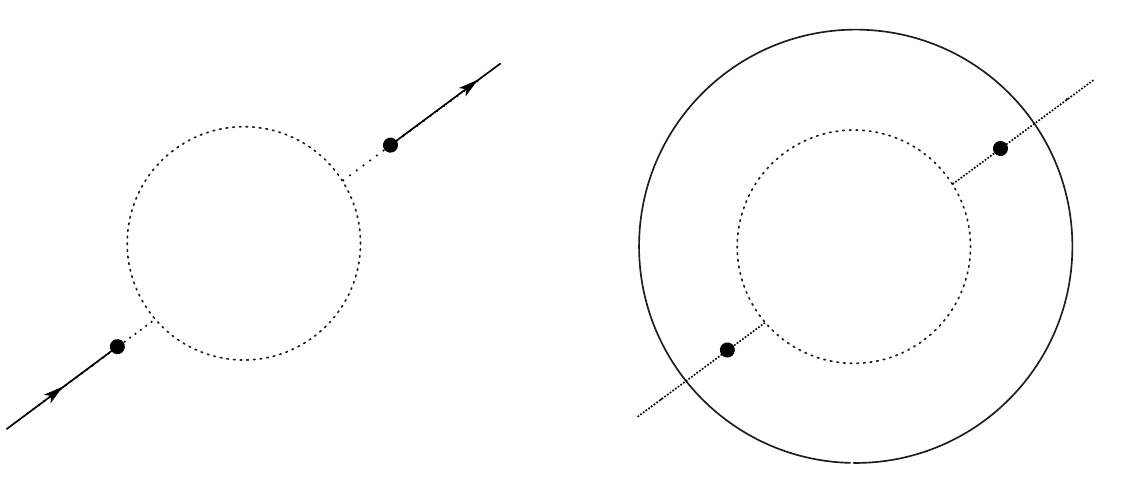
\caption{Contours $\Gm'$ and $\Gm''$, their orientations and $\pm$ sides.}
\label{disegno2}
\end{figure}

See Figure \ref{disegno2} for specification of the orientations of $\Gm'$ and $\Gm''$. We also denote by $\Om_\pm$ the half-planes defined by
\[\Om_-:=\{z\colon \tau-\pi<\arg z<\tau\},\quad \Om_+:=\{z\colon \tau<\arg z<\tau+\pi\},
\]and by $\Gm_0$ the oriented line $\Gm_0:=\Gm'\cup [P_1,P_2]$, where $[P_1,P_2]$ denotes the oriented segment from $P_1$ to $P_2$. We have $\der\Om_+=-\der\Om_-=\Gm_0$.
\vskip1,5mm
Second, following \cite[Ch.\,11]{Dur70}, we introduce some natural classes of analytic functions on which to work.
\begin{defn}
We define the {\it Smirnov class} $\mathscr E^p(\Om_{\pm})$, with $0<p<\infty$, as the set of analytic functions $f$ on $\Om_\pm$ such that
\[\sup_{\ell\in\frak R_\pm}\int_{\ell}|f(\zeta)|^p|d\zeta|<+\infty,
\]
where $\frak R_\pm$ is the set of lines $\ell\subseteq\Om_\pm$ parallel to $\Gm_0$.
\end{defn}
\begin{rem}
One can prove that if $f\in\mathscr E^p(\Om_+)$ (resp.\,\,$\Om_-$), with $0<p<\infty$, then $f(z)\to 0$ for $z\to\infty$ within each half-plane in $\Om_ +$ (resp.\,\,$\Om_-$) with boundary $\ell\in\frak R_+$ (resp.\,\,$\ell\in\frak R_-$). 
See \cite[Cor.\,2 of Th.\,11.3]{Dur70}.
\end{rem}
Essentially, the Smirnov classes $\mathscr E^p(\Om_\pm)$ are the largest classes of analytic functions (vanishing at $z=\infty$) for which Cauchy integral formula holds on $\Om_\pm$. More precisely, we have the following results.
\begin{thm}[{\cite[Th.\,11.1 and Th.\,11.8]{Dur70}}]\label{thdur}$\quad$
\begin{enumerate}
\item If $f\in\mathscr E^p(\Om_\pm)$, with $0<p<\infty$, then $f$ has non-tangential limit at the boundary $\Gm_0$ of $\Om_\pm$ for almost all $z_o\in\Gm_0$, and $f|_{\Gm_0}\in L^p(\Gm_0;|d\zeta|)$.
\vskip1,5mm
\item If $f\in\mathscr E^p(\Om_+)$ (resp. $\mathscr E^p(\Om_-)$), with $1\leq p<\infty$, then
\[f(z)=\frac{1}{2\pi\sqrt{-1}}\int_{\der\Om_+}\frac{f(\zeta)}{\zeta-z}d\zeta,\quad z\in\Om_+,\]
\[\left(\text{resp. }f(z)=\frac{1}{2\pi\sqrt{-1}}\int_{\der\Om_-}\frac{f(\zeta)}{\zeta-z}d\zeta,\quad z\in\Om_-,\right)
\]and the integral vanishes for $z\in \Om_-$ (resp. $\Om_+$).
\vskip1,5mm
\item Conversely, if $h\in L^p(\Gm_0;|d\zeta|)$, with $1\leq p<\infty$, and
\[\frac{1}{2\pi\sqrt{-1}}\int_{\der\Om_+}\frac{h(\zeta)}{\zeta-z}d\zeta\equiv 0,\quad z\in\Om_-, \quad (\text{resp. }z\in\Om_+,)
\]then for $z\in\Om_+$ (resp. $z\in\Om_-$) this integral represents a function $f\in\mathscr E^p(\Om_+)$ (resp. $\mathscr E^p(\Om_-)$) whose boundary function $f=h$ a.e.\qed
\end{enumerate}
\end{thm}

\begin{prob}[{Problem $\mc P_1[\bm u,\tau,\frak M]$}]Find an analytic function $\Psi\colon \C\setminus\Gm'\to M_n(\C)$ such that 
\begin{enumerate}
\item $(\Psi-I)|_{\Om_\pm}\in\mathscr E^2(\Om_\pm)$, so that the non-tangential limits $\Psi_\pm$ of $\Psi$ from the $-$ and $+$ sides of $\Gm'$ exist a.e., and $\Psi_{\pm}\in L^2(\Gm';|d\zeta|)$,
\item the limits $\Psi_\pm$ are related by
\beq
\label{rhba1}\Psi_+(z)=\Psi_-(z)H(z;\bm u).
\eeq

\end{enumerate}
\end{prob}

\begin{lem}\label{HL2}
We have $H(\zeta;\bm u)-I\to 0$ exponentially fast for $\zeta\to\infty$ along $\Gm'$. In particular, $H(-;\bm u)-I\in L^2(\Gm';|d\zeta|)$.
\end{lem}

\proof
The $(i,j)$-entry of $H(\zeta;\bm u)-I$ equals
\beq\label{entryij}
c_{ij}\exp\{(u^i-u^j)\zeta+(B_{ii}-B_{jj})\log\zeta\}-\delta_{ij},
\eeq
where 
 \begin{empheq}[left=c_{ij}{=}\empheqlbrace]{align*}
 &(S_1^{-1})_{ij},\quad\text{along }\ell_1,\\
 &(S_2)_{ij},\quad\text{along }\ell_2.
    \end{empheq}
By conditions $(5),(6),(7)$ (and $(9)$ if $\bm u\in\Dl$) of Definition \ref{datum}, we deduce that \eqref{entryij} goes to zero exponentially fast for $\zeta\to\infty$ along $\ell_1$ and $\ell_2$.
\endproof

\begin{thm}\label{solpa1}
If $\min\{|P_1|,|P_2|\}$ is sufficiently big, then there exists a unique solution $\Psi$ of the problem $\mc P_1[\bm u,\tau,\frak M]$, holomorphically depending on $\bm u\in\mc V$. Moreover, $\det \Psi\equiv 1$.
\end{thm}
\proof
If $\Psi$ is a solution of $\mc P_1[\bm u,\tau,\frak M]$, the condition $(\Psi-I)|_{\Om_\pm}\in\mathscr E^2(\Om_\pm)$ captures both the condition $\Psi(\infty)=I$, and the validity of a Cauchy integral representation. Indeed, we have
\begin{align*}
\Psi(z;\bm u)&=I+\int_{\Gm'}\frac{\Psi_+(\zeta)-\Psi_-(\zeta)}{\zeta-z}\frac{d\zeta}{2\pi\sqrt{-1}}\\
&=I+\int_{\Gm'}\frac{\Psi_-(\zeta)\left(H(\zeta;\bm u)-I\right)}{\zeta-z}\frac{d\zeta}{2\pi\sqrt{-1}},\quad z\notin \Gm',
\end{align*}by Theorem \ref{thdur}, and the jump condition \eqref{rhba1}. 
See also \cite[Ch. 3]{Pain}\cite[\S 5.1.3]{Its11} \cite[Ch. 2]{TO16}.
Set $\delta\Psi:=\Psi-I$ and $\delta H:=H-I$. The previous equation can be written as
\beq\label{eqint}
\delta\Psi(z)=\int_{\Gm'}\frac{\delta\Psi_-(\zeta)\delta H(\zeta;\bm u)}{\zeta-z}\frac{d\zeta}{2\pi\sqrt{-1}}+\int_{\Gm'}\frac{\delta H(\zeta;\bm u)}{\zeta-z}\frac{d\zeta}{2\pi\sqrt{-1}},\quad z\notin\Gm'.
\eeq
Given a function $f$ defined on $\Gm'$, introduce the functions $\mc C^\pm_{\Gm'}[f]$ on $\Gm'$ defined by the Cauchy integrals
\[\mc C^\pm_{\Gm'}[f](p):=\lim_{z\to p_{\pm}}\int_{\Gm'}\frac{f(\zeta)}{\zeta-z}\frac{d\zeta}{2\pi\sqrt{-1}},\quad p\in\Gm',
\]whenever the integral is finite. General results ensure that if $f\in L^p(\Gm';|d\zeta|)$, with $1\leq p<\infty$, then $\mc C^{\pm}_{\Gm'}[f]$ exists for $p\in\Gm'$ a.e.. Moreover, the Cauchy operators $\mc C^{\pm}_{\Gm'}$ are bounded in $L^p(\Gm';|d\zeta|)$, with $1<p<\infty$, i.e. there exists a constant $k_p>0$ such that
\[\|\mc C^{\pm}_{\Gm'}[f]\|_{L^p(\Gm')}\leq k_p\| f\|_{L^p(\Gm')},\quad \text{if }f\in L^p(\Gm';|d\zeta|).
\]
See \cite{MP80,Zho89,DZ02a,DZ02b,TO16} for details and proofs.
Taking the limit $z\to z_-$ in \eqref{eqint}, we obtain the following integral equation for $\dl\Psi_-$:
\beq
\label{eqint2}
\mc C[\dl H;\Gm']\dl\Psi_-=\mc C_{\Gm'}^-[\dl H],
\eeq
where
\[\mc C[\dl H;\Gm']f:=f-\mc C_{\Gm'}^-[f\cdot\dl H].
\]Notice that $\mc C_{\Gm'}^-[\dl H]\in L^2(\Gm';|d\zeta|)$ by Lemma \ref{HL2}. Moreover, if $f\in L^2(\Gm';|d\zeta|)$, we have
\[\|\mc C_{\Gm'}^-[f\cdot\dl H]\|_{L^2(\Gm')}\leq k_2\|f\cdot\dl H\|_{L^2(\Gm')}\leq k_2\cdot \sup_{\zeta\in\Gm'}\|\dl H(\zeta;\bm u)\|\cdot \|f\|_{L^2(\Gm')}.
\]By Lemma \ref{HL2}, we can assume that $\min\{|P_1|,|P_2|\}$ is so big that 
\[\sup_{\zeta\in\Gm'}\|\dl H(\zeta;\bm u)\|<\frac{1}{1+k_2},
\]then the operator $\mc C[\dl H;\Gm']\colon L^2(\Gm';|d\zeta|)\to L^2(\Gm';|d\zeta|)$ is invertible with inverse
\[\mc C[\dl H;\Gm']^{-1}=\sum_{m=0}^\infty \mc C_{\Gm'}^-[(-)\cdot\dl H]^m,
\]that is, for any $f\in L^2(\Gm';|d\zeta|)$, we have
\[\mc C[\dl H;\Gm']^{-1}f=f+\mc C_{\Gm'}^-[f\cdot\dl H]+\mc C_{\Gm'}^-[\mc C_{\Gm'}^-[f\cdot\dl H]\cdot\dl H]+\dots.\]Equation \eqref{eqint2} can be uniquely solved in $\dl\Psi_-$, and the formula \eqref{eqint} gives the unique solution $\Psi$ of the RHB boundary value problem. Notice that the Cauchy operator $\mc C^-_{\Gm'}[(-)\cdot\delta H]$ depends holomorphically on $\bm u$, so that $\Psi(z;\bm u)$ is holomorphic in $\bm u$.
Finally, notice that the jump condition \eqref{rhba1} implies
\[\det\Psi_+=\det\Psi_-\det H=\det\Psi_-,
\]since $\det H(\zeta;\bm u)\equiv 1$ along $\Gm'$. Hence, $\det\Psi$ is an entire function, and from the asymptotic condition $\Psi\to I$ for $z\to\infty$, we deduce $\det\Psi\equiv 1$ by Liouville Theorem.
\endproof
\begin{rem}
It follows from general results of \cite[Chapter 4]{Mus72} that the function $\Psi$ has two singularities at $P_1,P_2$. One can also prove that $\Psi_{\pm}\in L^2(\Gm';|d\zeta|)$ admit continuous representatives  on $\Gm'\setminus\{P_1,P_2\}$.
\end{rem}
\vskip3mm
Define the function $\frak S(-;\bm u)\colon \C\setminus\Gm\to GL(n,\C)$ by
 \begin{empheq}[left=\frak S(z;\bm u){:=}\empheqlbrace]{align*}
 &I,\quad\text{for }z\in\Pi_0,\\
 &H_{1}(z;\bm u)^{-1},\quad\text{for }z\in\Pi_R,\\
 &H_2(z;\bm u)^{-1},\quad\text{for }z\in\Pi_L. 
    \end{empheq}
    
    \begin{lem}
    The function $\frak S(-;\bm u)$ is a \emph{naive} solution of $\mc P[\bm u,\tau,\frak M]$: it satisfies conditions (1),(2), 
    but not  (3).
        \end{lem}
        \proof
        This is easily checked, by invoking equations \eqref{sc1}, \eqref{sc2}.
        \endproof
        
Consider the function $\widetilde H(-;\bm u)\colon \C\setminus \Gm\to GL(n,\C)$ defined by
\[\widetilde H(z;\bm u):=\Psi(z,\bm u)\frak S(z;\bm u)^{-1},
\]where $\Psi$ is the unique piecewise analytic solution of $\mc P_1[\bm u,\tau,\frak M]$, as in Theorem \ref{solpa1}.

\begin{lem}\label{lcht}$\quad$
\begin{enumerate}
\item The function $\widetilde H(-;\bm u)$ is continuous along $\Gm''$.
\item The function $\det\widetilde H(-;\bm u)$ has zero index across $\Gm''$, i.e.
\[{\rm ind}_{\Gm''}\det\widetilde H(-;\bm u):=\frac{1}{2\pi\sqrt{-1}}\oint_{\Gm''}d\log\det\widetilde H(\zeta;\bm u)=0.
\]
\end{enumerate}
\end{lem}
\proof
Point (1) is obvious. For point (2) notice that, for $\zeta\in\Gm''$, we have
\begin{align*}\log\det \widetilde H(\zeta;\bm u)&=\log\det \frak S(\zeta;\bm u)^{-1}\\
&=\log\det H_1(\zeta;\bm u)\\
&=\tr Q(\zeta;\bm u)-\log\det C-(\tr L+\tr D)\log \zeta\\
&=\zeta\sum_{i=1}^n u^i-\log\det C+\underbrace{(\tr B-\tr L-\tr D)}_{0\text{ by \eqref{vad1}}}\log \zeta.
\end{align*}
For the third line, see equation \eqref{defH}. This completes the proof.
\endproof

We can now introduce a second auxiliary RHB boundary value problem, with continuous coefficients on the simple closed contour $\Gm''$.

\begin{prob}[{Problem $\mc P_2[\bm u,\tau,\frak M]$}]Find an analytic function $\Upsilon\colon \C\setminus\Gm''\to M_n(\C)$ such that 
\begin{enumerate}
\item the non-tangential limits $\Upsilon_\pm\colon \Gm''\to M_n(\C)$ of $\Upsilon$ from the - and + sides of $\Gm''$ exist,
\item they are related by
\beq
\label{rhba2}\Upsilon_+(z)=\Upsilon_-(z)\widetilde H(z;\bm u),
\eeq
\item $\Upsilon(z)$ tends to the identity matrix $I$ as $z\to\infty$.
\end{enumerate}
\end{prob}

\begin{thm}\label{teosolv1}
The solvability of $\mc P[\bm u,\tau,\frak M]$ is equivalent to the solvability of $\mc P_2[\bm u,\tau,\frak M]$. 
\end{thm}
\proof
If $G$ is the solution of $\mc P[\bm u,\tau,\frak M]$, then 
 \begin{empheq}[left=\Upsilon (z;\bm u){:=}\empheqlbrace]{align*}
 &G(z;\bm u)\Psi(z;\bm u)^{-1},\quad \text{for $z$ outside $\Gm''$,}\\
 &G(z;\bm u)\frak S(z;\bm u)^{-1},\quad \text{for $z$ inside $\Gm''$,}
    \end{empheq}
    is the solution of $\mc P_2[\bm u,\tau,\frak M]$. Vice-versa, if $\Upsilon$ is the solution of $\mc P_2[\bm u,\tau,\frak M]$, then the solution $G$ of $\mc P[\bm u,\tau,\frak M]$ is obtained by inverting the equations above.
\endproof

\subsection{Solvability as an open property}If $\Upsilon$ is a solution of $\mc P_2[\bm u,\tau,\frak M]$, then we have
\beq
\label{if2}
\Upsilon(z)=I+\int_{\Gm''}\frac{\Upsilon_-(\zeta)(\widetilde H(\zeta;\bm u)-I)}{\zeta-z}\frac{d\zeta}{2\pi\sqrt{-1}}.
\eeq
In the limit $z\mapsto z_-$, we obtain the integral equation
\beq
\label{ie2}
\Upsilon_-=I+\mc C^-_{\Gm''}\left[\Upsilon_-\delta\widetilde H\right],\quad \delta\widetilde{H}:=\widetilde H-I,
\eeq
where $\mc C^\pm_{\Gm''}$ denotes the Cauchy integrals with respect to the contour $\Gm''$.
Conversely, if $\Upsilon_-$ is a solution of \eqref{ie2}, then \eqref{if2} gives the solution of $\mc P_2[\bm u,\tau,\frak M]$, see \cite[Ch. 3]{Pain}\cite[\S 5.1.3]{Its11}\cite[Ch. 2]{TO16}. 

\begin{thm}
The operator 
\[T(\bm u)\colon L^2(\Gm'';|d\zeta|)\to L^2(\Gm'';|d\zeta|),\quad f\mapsto f-\mc C^-_{\Gm''}[f\cdot \delta\widetilde H(\bm u)]
\]is a Fredholm operator with index\footnote{Recall that the index of a Fredholm operator $T$ is the integer ${\rm ind}\,T:=\dim\ker T-\dim {\rm coker}\,T$.} $0$.
\end{thm}

\proof
Assume we are given a factorization $\widetilde H(\zeta;\bm u)=(I-W^-(\zeta;\bm u))^{-1}(I+W^+(\zeta;\bm u))$ with 
\bea
W^+(z;\bm u)\in L^\infty(\Gm'')\cap L^2(\Gm''),\qquad (I- W^-(\zeta;\bm u))^{- 1}-I\in L^\infty(\Gm'')\cap L^2(\Gm'').
\eea
Define the Cauchy type operator 
\[\mc C_W\colon L^2(\Gm'')\to L^2(\Gm''),\quad f\mapsto \mc C^+_{\Gm''}[fW^-]+\mc C^-_{\Gm''}[fW^+].
\]Standard results imply that the operator $f\mapsto f-\mc C_W[f]$ is Fredholm, and its index is given by
\[{\rm ind}(Id-\mc C_W)={n}\,{\rm ind}_{\Gm''}\det\widetilde H=0,
\]by Lemma \ref{lcht} point (2), see \cite{Zho89,TO16}. In our case, we can take $W^-=0$ and $W^+=\dl\widetilde H$ by Lemma \ref{lcht} point (1). This completes the proof.
\endproof

A set $\Theta\subseteq \C^n$ is said to be {\it analytic} if it is the zero locus of a scalar analytic function.

\begin{thm}\label{teoimp}
Let $\bm u_o\in\C^n$. Assume that the pair $(\tau,\frak M)$ is admissible at each point of a sufficiently small open neighborhood $\mc V$ of $\bm u_o$. If $\mc P[\bm u_o,\tau,\frak M]$ is solvable, there exists an analytic set $\Theta\subseteq \mc V\setminus \{\bm u_o\}$ such that $\mc P[\bm u,\tau,\frak M]$ is solvable for all $\bm u\in\mc V\setminus\Theta$. Moreover, the solution $G(z;\bm u)$ is holomorphic with respect to $\bm u\in \mc V\setminus\Theta$.
\end{thm}

For the proof we firstly invoke the following Lemma. 

\begin{lem}[{\cite[Lemma 10]{Gra70}}]\label{lemgra} Let $X$ be a Banach space and $\frak F(X)$ be the set of its Fredholm operators. Let $\Om\subseteq\C^n$ be a connected domain, and $T\colon\Om\to \frak F(X)$ a holomorphic function. If $T(\bm u_o)^{-1}$ exists for some $\bm u_o\in\Om$, then $T(\bm u)^{-1}$ exists on the complement $\Om\setminus \Theta$ of an analytic set $\Theta$, and $T^{-1}$ is meromorphic on $\Om$.\qed
\end{lem}

\begin{rem}
Lemma \ref{lemgra} was originally due to I.\,Gohberg and E.\,Sigal in the case $n=1$, \cite{GS70}. The general case was proved by B.\,Gramsch, though special cases were previously obtained by several authors. For a sketch of a proof, based on arguments of \cite{GS70} and \cite[XI.8]{GGK90},  see \cite[Sec. 2]{Kab12}.
\end{rem}

\proof[Proof of Theorem \ref{teoimp}]
By assumption and Theorem \ref{teosolv1}, the $\mc P_2[\bm u_o,\tau,\frak M]$ is solvable. We claim that the solution $\Upsilon$ is unique. The function $\det\Upsilon(z;\bm u_o)$ solves the scalar RH problem 
\[\det\Upsilon_+(z;\bm u_o)=\det\Upsilon_-(z;\bm u_o)\det\widetilde H(z;\bm u_o).
\]Since the function $\det\widetilde H(-;\bm u_o)$ has zero index along $\Gm''$, this scalar equation can be uniquely solved: the solution is given by
\[\det\Upsilon(z;\bm u_o)=\exp\int_{\Gm ''}\frac{\log\det \widetilde H(\mu;\bm u_o)}{\mu-z}\frac{d\mu}{2\pi\sqrt{-1}},
\]see e.g. \cite[\S 2.3.1]{TO16}. In particular, $\Upsilon(z;\bm u_o)$ is invertible. Assume that $\Upsilon(z;\bm u_o),\tilde\Upsilon(z;\bm u_o)$ are two solutions of $\mc P_2[\bm u_o,\tau,\frak M]$. Put $X(z):=\Upsilon(z;\bm u_o)\tilde\Upsilon(z;\bm u_o)^{-1}$. For $z\in\Gm''$ we have
\[X_+(z)=\Upsilon_+(z;\bm u_o)\tilde\Upsilon_+(z;\bm u_o)^{-1}=\Upsilon_-(z;\bm u_o)\widetilde H(z;\bm u_o)\widetilde H(z;\bm u_o)^{-1}\tilde\Upsilon_-(z;\bm u_o)^{-1}=X_-(z).
\]Hence $X(z)$ is analytic, and moreover $X(z)\to I$ for $z\to\infty$. By Liouville Theorem we have $X(z)\equiv I$, and $\Upsilon=\tilde \Upsilon$.
It follows that the Fredholm operator $T(\bm u_o)$ has both trivial kernel and index zero. Hence $T(\bm u_o)^{-1}$ exists, Lemma \ref{lemgra} applies, and the problem $\mc P_2[\bm u,\tau,\frak M]$ is solvable on the complement of an analytic set $\Theta\subseteq\mc V\setminus \{\bm u_o\}$. By Theorem \ref{teosolv1} one concludes.
\endproof

\subsection{Proof of Sabbah Theorem}\label{ssabt}Let $(E_o,\nabla_o)$ be in Birkhoff normal form \eqref{con1} with $\La_o={\rm diag}(u_o^1,\dots,u_o^n)$ and $\bm u_o\in\Dl$. Consider the differential system defining $\nabla_o$-flat sections 
\beq
\label{sist1}
\frac{dY}{dz}=\left(\La_o+\frac{1}{z}B_o\right)Y,
\eeq where $Y$ is a matrix-valued function. 

\begin{prop}[\cite{AB}{\cite{CG18}\cite[Section 16]{CDG1}}]\label{sol0}
The differential system \eqref{sist1} has a fundamental system of solutions in \emph{Birkhoff-Levelt normal form}
\[Y_0(z)=\mc G_0(z)z^{\mc D}z^{S+R},\quad \mc G_0(z)=K\left(I+\sum_{j=1}^\infty A_jz^j\right),
\]where
\begin{itemize}
\item $K$ puts $B_o$ in Jordan form $J=K^{-1}B_oK$,
\item $\mc D$ is a diagonal matrix of integers (called \emph{valuations}),
\item $S$ is a Jordan matrix whose eigenvalues have real part in $[0,1[$,
\item $R$ is a nilpotent matrix, with possibly non-vanishing entries only if some of the eigenvalues of the matrix $B_o$ differ by a non-zero integer.
\end{itemize}
Moreover, we have
\[\pushQED{\qed} 
J={\mc D}+S.\qedhere
\popQED
\]
\end{prop}

\begin{prop}[{\cite[Prop. 4.2]{CDG1}}]\label{solinf}
Assume that 
\begin{enumerate}
\item $B''_o\in {\rm Im\ ad}(\La(\bm u_o))$,
\item $B'_o$ is partially non-resonant, i.e. $(B_o)_{ii}-(B_o)_{jj}\notin\Z\setminus\{0\}$ if $u_o^i= u_o^j$.
\end{enumerate}
Then, the differential system \eqref{sist1} has a unique formal solution of the form
\[Y_F(z)=\left(I+\sum_{k=1}^\infty F_k z^{-k}\right)z^{B_o'}e^{\La_oz}.
\]
If $\tau$ is admissible at $\bm u_o$, then there exist three fundamental systems of solutions $Y_1,Y_2,Y_3$ of \eqref{sist1}, satisfying respectively
\beq\label{asym}
Y_h(z)\sim Y_F(z),\quad |z|\to+\infty,\quad \tau-(3-h)\pi<\arg z<\tau+(h-2)\pi,\quad h=1,2,3.
\eeq
and uniquely determined by these conditions.
\qed
\end{prop}
\begin{rem}
For $h=1,2,3$, set
\[\mc G_h(z):=Y_h(z)e^{-\La_oz}z^{-B_o'},
\]
\[\mc W_{\tau,h}:=\left\{z\in\widetilde{\mathbb C^*}\colon \tau-(3-h)\pi<\arg z<\tau+(h-2)\pi\right\}.
\]Denote by $\overline{\mathcal W}$ an arbitrary unbounded closed sector, in the universal cover $\widetilde{\mathbb C^*}$ of $\C^*$, with vertex at $\infty$.
The precise meaning of the asymptotic relations \eqref{asym} is the following: \[
\forall h\in\{1,2,3\},
\
 \forall \ell\in\mathbb N,
 \
  \forall \overline{\mathcal W}\subsetneq \mathcal W_{\tau,h},
  \
   \exists C_{h,\ell,\overline{\mathcal W}}>0\colon \text{ if }z\in\overline{\mathcal W}\setminus\left\{0 \right\}\text{ then }\]
  \[
    \left\| \mc G_h(z)-\left(I+\sum_{m=1}^{\ell-1}\frac{F_m}{z^m}\right)\right\|
    <\frac{C_{h,\ell,\overline{\mathcal W}}}{|z|^\ell}.
\]
\end{rem}
\vskip3mm
In the notations of Propositions \ref{sol0}, \ref{solinf}, consider the $6$-tuples $\frak M=(B,D,L,S_1,S_2,C)$ where
\[B:=B_o',\quad D:=\mc D,\quad L:= S+R,\quad 
\]
and the matrices $S_1,S_2, C$ are defined by 
\beq
\label{18.05.21-1}
Y_2(z)=Y_1(z) S_1,\quad Y_3(z)= Y_2(z)S_2,\quad Y_2(z)=Y_0(z)C.
\eeq

\begin{prop}\label{18.05.21-2}
The $6$-tuple $\frak M$ is a $(\bm u_o,\tau)$-admissible datum. The RHB boundary value problem $\mc P[\bm u_o,\tau,\frak M]$ is solvable, with unique solution 
 \begin{empheq}[left=G(z;\bm u_o){=}\empheqlbrace]{align*}
 &\mc G_0(z),\quad z\in\Pi_0,\\
 &\mc G_2(z),\quad z\in\Pi_R,\\
 &\mc G_3(z),\quad z\in\Pi_L.
   \end{empheq}
\end{prop}
\proof
Conditions (1),(2),(3) of Definition \ref{datum} are trivially satisfied. The proof of conditions (4),(5),(6),(7) for the Stokes matrices $S_1,S_2$ are standard, see e.g. \cite[Section 6.3]{CDG1}. Denote by $\widetilde{\mathbb C^*}$ the universal cover of $\C^*$. Notice that 
\[Y_3(ze^{2\sqrt{-1}\pi})=Y_1(z)e^{2\sqrt{-1}\pi B_o'}, \quad z\in\widetilde{\C^*},
\]both sides having the same asymptotic expansion $Y_F(ze^{2\sqrt{-1}\pi})$ for $|z|\to+\infty$, and $\tau-2\pi<\arg z<\tau-\pi$. We deduce that
\[Y_0(ze^{2\sqrt{-1}\pi})CS_2=Y_0(z)CS^{-1}_1e^{2\sqrt{-1}\pi B_o'},
\]so that
\[e^{2\sqrt{-1}\pi L}=CS^{-1}_1e^{2\sqrt{-1}\pi B_o'}S_2^{-1}C^{-1}.
\]This proves condition (8) of Definition \ref{datum}. Finally, condition (9) follows from \cite[Th. 2.1]{CG18}, \cite[Prop. 6.1]{CDG1}. The remaining part of the statement follows from equations \eqref{18.05.21-1}, and the uniqueness stated in Proposition \ref{solinf}.
\endproof

By Proposition \ref{18.05.21-2} and Theorem \ref{teoimp}, there exist an open neighborhood $\mc V$ of $\bm u_o$, an analytic set $\Theta\subseteq\mc V\setminus\{\bm u_o\}$ on which the RHB problem $\mc P[\bm u,\tau,\frak M]$ is solvable, with unique solution $G(z;\bm u)$ holomorphic with respect to $\bm u\in\mc V\setminus\Theta$. Define the functions
\bea
&Y_{2}(z;\bm u):=G(z;\bm u)z^{B'_o}e^{\La(\bm u)z},\quad &z\in\Pi_{R},\\
&Y_{3}(z;\bm u):=G(z;\bm u)z^{B'_o}e^{\La(\bm u)z},\quad &z\in\Pi_{L},\\
&Y_{0}(z;\bm u):=G(z;\bm u)z^{\mc D}z^L,\quad &z\in\Pi_{0}.
\eea
By Remark \ref{26.05.21-rem1}, we have
\[
Y_2(z;\bm u)=Y_0(z;\bm u)\cdot C,\quad Y_3(z;\bm u)=Y_2(z;\bm u)\cdot S_2.
\]
Moreover, we have $G(z;\bm u)=I+\frac{F_1(\bm u)}{z}+O\left(\frac{1}{z^2}\right)$ in $z\to\infty$ in $\Pi_{L/R}$, so that
\bea
\frac{\der Y_{2/3} }{\der z}\cdot Y_{2/3}^{-1}&=&\der_zG\cdot G^{-1}+\frac{1}{z}GB_o'G^{-1}+G\La G^{-1}\\
&=&\La(\bm u)+\frac{1}{z}\left([F_1(\bm u),\La(\bm u)]+B'_o\right)+O\left(\frac{1}{z^2}\right),\quad z\to\infty,\\
\frac{\der Y_{0} }{\der z}\cdot Y_{0}^{-1}&=&\der_zG\cdot G^{-1}+\frac{1}{z}\left(GDG^{-1}+Gz^DLz^{-D}G^{-1}\right)\\
&=& \frac{1}{z}K\left(D+L\right)K^{-1}+O(1),\quad z\to0.
\eea
The matrices $S_1,S_2,C$ are constant with respect to both $\bm u$ and $z$: we deduce that the r.h.s.\,\,of the two equalities above are equal. This implies that $Y_{2}, Y_3$, and $Y_0$ are solutions of the differential equation
\beq
\label{mde1}
\frac{\der}{\der z}Y=\left[\La(\bm u)+\frac{1}{z}V(\bm u)\right]Y,\quad V(\bm u):=[F_1(\bm u),\La(\bm u)]+B'_o.
\eeq
Similarly, we have
\bea
\frac{\der Y_{2/3} }{\der u^i}\cdot Y_{2/3}^{-1}&=&\frac{\der G}{\der u^i}\cdot G^{-1}+zGE_iG^{-1}=zE_i+[F_1,E_i]+O\left(\frac{1}{z}\right),\\
\frac{\der Y_{0} }{\der u^i}\cdot Y_{0}^{-1}&=&\frac{\der G}{\der u^i}\cdot G^{-1}=\frac{\der G_0}{\der u^i}\cdot G^{-1}_0+O(z),
\eea
where $(E_i)_{ab}=\dl_{ai}\dl_{bi}$ and $G(z;\bm u)=G_0(\bm u)+O(z)$ for $z\to 0$ (and in particular $G_0(\bm u_o)=K$). The matrices $S_1,S_2,C$ being constant, we deduce that the r.h.s.\,\,of the two equalities above are equal. Hence $Y_2,Y_3$, and $Y_0$ are solutions of the differential systems
\beq
\label{mde2}
\frac{\der}{\der u^i}Y=\left(zE_i+V_i(\bm u)\right)Y,\quad V_i(\bm u):=[F_1(\bm u),E_i]=\frac{\der G_0}{\der u^i}\cdot G^{-1}_0,\quad i=1,\dots,n.
\eeq
The datum of the compatible joint differential systems \eqref{mde1} and \eqref{mde2}, for $\bm u\in\mc V\setminus \Theta$, proves the statement of Sabbah Theorem \ref{tsab}. 

\begin{rem}
Note that in equations \eqref{mde1} and \eqref{mde2} we can replace $F_1$ with its off-diagonal part $F_1''$, since both $\La(\bm u)$ and $E_i$ are diagonal.
\end{rem}

\begin{rem}
Propositions \ref{sol0} and \ref{solinf} also hold true for $\bm u_o\in\C^n\setminus\Dl$, these are standard results. All the subsequent arguments can be applied, giving an analytical proof of Theorem \ref{tmal}.
\end{rem}

\section{Formal Frobenius and \dfms}\label{sec4}
We briefly review basic notions of the theory of Frobenius manifolds, in both formal and analytic frameworks. General references are \cite{Dub96,Dub98,Dub99,man,Her02,Sab07}.

\subsection{Formal Frobenius manifolds}\label{ffm}
Let 
\begin{itemize}
\item $k$ be a commutative $\Q$-algebra, 
\item $H$ be a free $k$-module of finite rank, 
\item $\eta\colon H\otimes H\to k$ be a symmetric pairing, inducing an isomorphism $\eta'\colon H\to H^T$, where $H^T$ is the dual module,
\item $K:=k[\![H^T]\!]$ be the completed symmetric algebra of $H^T$.
\end{itemize}
Fix a basis $(\Dl_1,\dots, \Dl_n)$ of $H$, and denote by $\bm t=(t^1,\dots, t^n)$ the dual coordinates. The algebra $K$ is then identified with the algebra of formal power series $k[\![\bm t]\!]$. Denote by ${\rm Der}_k(K)$ the $K$-module of $k$-linear derivations of $K$. Put $\der_\al=\frac{\der}{\der t^\al}\colon K\to K$. It is well known that ${\rm Der}_k(K)$ is a free $K$-module with basis $(\der_1,\dots,\der_n)$, see e.g.\,\,\cite{Now86}. If $\Phi\in K$, we will write $\Phi_\al$ for $\der_\al\Phi$.
\vskip 2mm
Elements of $H_K:=K\otimes_k H$ will be identified with derivations on $K$, by $\Dl_\al\mapsto \der_\al$. 
\vskip2mm
For $\al,\bt=1,\dots,n$, set $\eta_{\al\bt}:=\eta(\Dl_\al,\Dl_\bt)$. The matrix $(\eta^{\al\bt})$ will denote the inverse of the Gram matrix $(\eta_{\al\bt})$ of $\eta$. Einstein summation rule will be used over repeated Greek indices.
\begin{defn}
A \emph{formal Frobenius manifold structure} on $(H,\eta)$ is given by a formal power series  $\Phi\in K$, called \emph{WDVV potential}, such that
\beq\label{wdvv}
\Phi_{\al\bt\gm}\eta^{\gm\dl}\Phi_{\dl\eps\phi}=\Phi_{\phi\bt\gm}\eta^{\gm\dl}\Phi_{\dl\eps\al},\quad \al,\bt,\eps,\phi=1,\dots,n.
\eeq
\end{defn}

Define the $K$-linear multiplication $\circ$ on $H_K$ by
\beq\label{fprod}
\Dl_\al\circ \Dl_\bt:=c_{\al\bt}^{\gm}\Dl_\gm,\quad \al,\bt=1,\dots, n,
\eeq
where $c_{\al\bt}^{\gm}:=\Phi_{\al\bt\dl}\eta^{\dl\gm}$. The WDVV equations \eqref{wdvv} are equivalent to the associativity of $\circ$.
\vskip2mm
An element $e\in H_K$ is called \emph{identity} if it is the identity for $\circ$. It is called \emph{flat identity} if $e\in H$. 
An element $E\in H_K$ is called \emph{Euler} if 
\bean
\label{eu1}
&&\frak L_E\eta=D\eta,\quad D\in k,\\
\label{eu2}
&&\frak L_E c= c.
\eean
Here $\eta$ is $K$-bilinearly extended to $H_K$. Moreover, $\frak L_E$ denotes the Lie derivative along $E$, and it is extended to the full tensor algebra of the $K$-module $H_K$, as follows\footnote{For this standard algebraic approach to Lie derivatives of tensors, see e.g.\,\,\cite{AMR88}}. Set 
\[{\rm T^0_0}=K,\quad\quad {\rm T}^{p}_q:=H_K^{\otimes p}\otimes_K \Hom_K(H_K, K)^{\otimes q},\quad p,q\geq 1,\quad \quad {\rm T}:=\bigoplus_{p,q\geq 0}{\rm T}^p_q.
\]We extend $\frak L_E$ on $\rm T$, by requiring that:
\begin{enumerate}
\item $\frak L_E$ acts on elements $f\in {\rm T^0_0}=K$ by $\frak L_Ef=Ef$;
\item $\frak L_E$ acts on elements $X\in {\rm T^1_0}=H_K\cong {\rm Der}_k(K)$ by $\frak L_EX=[E,X]$;
\item $\frak L_E$ is a tensorial derivation, i.e. it satisfies the Leibniz rule
\[\frak L_E(S_1\otimes S_2)=(\frak L_ES_1)\otimes S_2+S_1\otimes(\frak L_ES_2),\quad S_1\in{\rm T}^{p_1}_{q_1},\,S_2\in{\rm T}^{p_2}_{q_2}; \]
\item $\frak L_E$ commutes with contractions, i.e.
\bea
&&\frak L_E[V(\al_1,\dots,\al_r,X_1,\dots, X_s)]=(\frak L_EV)(\al_1,\dots,\al_r,X_1,\dots, X_s)\\
&&+\sum_{j=1}^r V(\al_1,\dots,\frak L_E \al_j,\dots,\al_r, X_1,\dots, X_s)
+\sum_{k=1}^s V(\al_1,\dots,\al_r, X_1,\dots, \frak L_EX_k,\dots, X_s),
\eea
for any $V\in {\rm T}^r_s$, $\al_1,\dots, \al_r\in{\rm T}^0_1$, $X_1,\dots, X_s\in{\rm T}^1_0$.
\end{enumerate}
In this paper we always consider Frobenius manifolds equipped with a flat identity $e=\Dl_1$, and an Euler element $E$. Equations \eqref{eu1} and \eqref{eu2} are then equivalent to the single equation
\[E\Phi=(1+D)\Phi+\text{quadratic terms in }\bm t.
\]

\begin{lem}\label{affeu}
Let $E=\sum_\al E^\al(\bm t)\Dl_\al$ be an Euler element. The series $E^\al(\bm t)\in K$ are linear affine expressions in $\bm t$, i.e. $E^\al(\bm t)=\sum_\bt A^\al_\bt t^\bt+r^\al$, for suitable coefficients $A^\al_\bt, r^\al\in k$.
\end{lem}
\proof
The Killing-conformal condition \eqref{eu1} reads $\eta_{\al\bt}\der_\gm E^\al+\eta_{\gm\al}\der_\bt E^\al= D\eta_{\bt\gm}$. Differentiating this equation by $\der_\la$, and permuting indices, we obtain
\[\der_\la\der_\gm E_\bt+\der_\la\der_\bt E_\gm=0, \quad \der_\gm\der_\bt E_\la+\der_\gm\der_\la E_\bt=0,\quad \der_\bt\der_\la E_\gm+\der_\bt\der_\gm E_\la=0,
\]where $E_\la:=\sum_\nu E^\nu\eta_{\nu\la}$. Since $\der_i\der_j=\der_j\der_i$, we conclude that $\der^2_{\la\bt}E_\gm=\der^2_{\la\gm}E_\bt=\der^2_{\bt\gm}E_\la=0$. It follows that ($E_\gm$, and consequently) $E^\gm$ are linear functions in $\bm t$.
\endproof

\subsection{\dfms} Given a complex analytic manifold $M$, we denote by $TM,T^*M$ its holomorphic tangent and cotangent bundles. If $E$ is a vector bundle on $M$, its $k$-th symmetric power is denoted by $\bigodot^kE$. 
\vskip2mm
A \emph{\dfm} structure on a complex manifold $M$ of dimension $n$ is defined by giving
\begin{enumerate}
\item[(FM1)] a symmetric non-degenerate $\mathcal O(M)$-bilinear form $\eta\in\Gamma\left(\bigodot^2T^*M\right)$, called {\it metric}\footnote{In what follows, the musical isomorphisms with respect to the metric $\eta$ will be denoted by $(-)^\flat$ and $(-)^\sharp$, respectively. If $\xi\in \Gm(TM)$, the 1-form $\xi^\flat\in\Gm(T^*M)$ is defined by $\xi^\flat(X)=\eta(X,\xi)$, where $X\in\Gm(TM)$. Conversely, if $\xi\in\Gm(T^*M)$, the vector field $\xi^\sharp\in\Gm(TM)$ is uniquely defined by the identity $\xi(X)=\eta(X,\xi^\sharp)$, where $X\in\Gm(TM)$. Thus $(-)^\flat\colon \Gm(TM)\to \Gm(T^*M)$ and $(-)^\sharp\colon \Gm(T^*M)\to \Gm(TM)$ are mutually inverse. In components, these operations are also known as ``lowering'' and ``raising'' indices, respectively. These operations naturally extend to mixed tensors.}, whose corresponding Levi-Civita connection $\nabla$ is flat;
\item[(FM2)] a $(1,2)$-tensor $c\in\Gamma\left(TM\otimes\bigodot^2T^*M\right)$ such that
\begin{enumerate}
\item the induced multiplication of vector fields $X\circ Y:=c(-,X,Y)$, for $X,Y\in\Gamma(TM)$, is \emph{associative}, 
\item $c^\flat\in\Gamma\left(\bigodot^3T^*M\right)$,
\item $\nabla c^\flat\in\Gamma\left(\bigodot^4T^*M\right)$;
\end{enumerate}
\item[(FM3)] a vector field $e\in\Gamma(TM)$, called the \emph{unity vector field}, such that
\begin{enumerate}
\item the bundle morphism $c(-,e,-)\colon TM\to TM$ is the identity morphism,
\item $\nabla e=0$;
\end{enumerate}
\item[(FM4)] a vector field $E\in\Gamma(TM)$, called the \emph{Euler vector field}, such that
\begin{enumerate}
\item $\frak L_Ec=c$,
\item $\frak L_E\eta=(2-d)\cdot \eta$, where $d\in\mathbb C$ is called the \emph{charge} of the Frobenius manifold.
\end{enumerate}
\end{enumerate}
\dfm s will be also called \emph{analytic} Frobenius manifolds.
\vskip2mm
By axiom (FM1), there exist systems of \emph{flat coordinates} $\bm t=(t^1,\dots, t^n)$, with respect to which the Levi-Civita connection $\nabla$ coincides with partial derivatives $\der_\al:=\frac{\der}{\der t^\al}$, for $\al=1,\dots,n$. Without loss of generality, we assume that the coordinate $t^1$ is such that $\der_1=e$. 
\vskip2mm
A \emph{pointed \dfm} is a pair $(M,p)$, where $M$ is a \dfm, and $p\in M$ is a fixed \emph{base point}. Given $(M,p)$ we will always consider flat coordinates $\bm t=(t^1,\dots,t^n)$ vanishing at $p$.

\subsection{From Dubrovin-Frobenius to formal Frobenius structures, and vice-versa}\label{affa} Given a pointed \dfm\ $(M,p)$, we can associate to it a formal Frobenius structure $(H,\eta,\Phi)$ over $k=\C$. Choose flat coordinates $\bm t$ vanishing at $p$, and set $H:=T_pM$ equipped with the metric $\eta|_p$. 
By axiom (FM2-c), the tensor $\der_\al c_{\bt\gm\dl}$ is completely symmetric: hence we deduce the local existence of a function $F$ such that $\der^3_{\al\bt\gm}F=c_{\al\bt\gm}$. By axioms (FM2-a), (FM2-b), we deduce that $F$ is a solution of WDVV equations, i.e.
\[\der_{\al\bt\gm}^3F\ \eta^{\gm\dl}\ \der_{\dl\eps\phi}^3F=\der_{\phi\bt\gm}^3F\ \eta^{\gm\dl}\ \der_{\dl\eps\al}^3F,\quad \al,\bt,\eps,\phi=1,\dots,n.
\]
Moreover, axiom (FM4) can be rephrased as the single equation
\[EF=(3-d)F+\text{quadratic terms in }\bm t.
\]
Let $\mc O_{M,p}$ be the local ring of germs at $p$, and $\frak m$ be its maximal ideal. The formal potential $\Phi$ is given by the image of $F$ in the completion $\widehat{\mc O_{M,p}}:=\varprojlim \left(\mc O_{M,p}/{\frak m}^\ell\right)$ of the local ring $\mc O_{M,p}$: this means that $\Phi$ is defined by the Taylor series expansion of $F$ at $p$ in coordinates $\bm t$. Moreover, the formal Frobenius structure $(H,\eta,\Phi)$ is also equipped with a flat unit $e|_p$ and an Euler vector field $E|_p$. We will say that the formal Frobenius structure constructed in this way, starting from a pointed \dfm, is \emph{convergent}.

Conversely, let us assume that $(H,\eta,\Phi)$ is a formal Frobenius structure over $k=\C$, with flat identity $e=\Dl_1$, and Euler element $E$. If the domain of convergence $\Om\subseteq H$ of the power series $\Phi\in k[\![\bm t]\!]$ is non-empty, then $\Om$ is equipped with a \dfm\  structure. The product $\circ$ defined in equation \eqref{fprod}, indeed, turns out to have analytic structure constants $c^\gm_{\al\bt}=\Phi_{\al\bt\la}\eta^{\la\gm}$. These are the components of a well-defined holomorphic section $c$ of the bundle $T\Om\otimes \bigodot^2T^*\Om$ as in axioms (FM2). The metric $\eta$ and the potential $\Phi$ are related by $\eta_{\al\bt}=\der^3_{1\al\bt}\Phi$. The unit element $e\in H_K$ is identified with the unit vector field $\der_1\in\Gm(T\Om)$. The Euler element $E\in H_K$ has linear components, by Lemma \ref{affeu}, so it can be identified with a well-defined holomorphic vector field on $\Om$, satisfying axioms (FM4).

\subsection{Semisimplicity of Frobenius structures.} In this Section we collect main results and properties which hold true for a wide class of Frobenius structures (both formal and analytic), namely \emph{semisimple} Frobenius structures.  We begin our exposition with the formal case.

\vskip3mm
Let $(H,\eta,\Phi)$ be a formal Frobenius manifold, and denote by $\circ_0$ the product on $H$ with structure constants $\Phi_{\al,\bt}^\gm(0)$. We say that $(H,\eta,\Phi)$ is 
\begin{itemize}
\item \emph{semisimple at the origin} if the $k$-algebra  $(H,\circ_0)$  is isomorphic to $k^n$;
\item \emph{formally semisimple} if the $K$-algebra $(H_K,\circ)$ is isomorphic to $K^n$.
\end{itemize}
In the first (resp. second) case there exist an idempotent basis $(\pi_1,\dots,\pi_n)$ of $H$ (resp. $H_K$) such that
\beq\label{idemp}
\pi_i\circ \pi_j=\pi_i\delta_{ij},\quad
 \eta(\pi_i,\pi_j)=0,\quad i\neq j.	
 \eeq
 Notice that the idempotent vectors $\pi_i$ are uniquely defined up to re-ordering.
\vskip2mm
\begin{lem}
A formal Frobenius manifold $(H,\eta,\Phi)$ is formally semisimple if and only if it is semisimple at the origin.
\end{lem}
\proof
Formal semisimplicity clearly implies semisimplicity at the origin. Let us prove the converse.  Denote by $\frak m:=(t^1,\dots,t^n)$ the maximal ideal of $K$. We will denote by $O(\frak m^p)$ an arbitrary sum of elements of $\frak m^p\cdot H_K$. 
For any fixed $h\in\N$ we call an \emph{$h$-order idempotent  basis} of $H_K$ a basis $(\pi_1^h,\dots,\pi_n^h)$ such that
\[\pi_i^h\circ\pi_i^h=\pi_i^h+O(\frak m^{h+1}),\quad\quad 
\pi_i^h\circ\pi_j^h=O(\frak m^{h+1}),
\]for $i,j=1,\dots,n$ and $i\neq j$.
Assume that $(H,\eta,\Phi)$ is semisimple at the origin. We claim there exist a $h$-order idempotent  basis of $H_K$ for any $h\in\N$. We prove it by induction on $h$. For $h=0$, it is trivial: if $(\pi_1^0,\dots,\pi_n^0)$ is an idempotent basis of $(H,\circ_0)$, then it is a 0-order idempotent basis of $H_K$.
Assume that $(\pi_1^h,\dots, \pi_n^h)$ is an $h$-order idempotent basis of $H_K$: we have
\[\pi_i^h\circ\pi_i^h=\pi_i^h+\sum_{k}a_{ik}\pi_k^h,\quad a_{ij}\in\frak m^{h+1},\quad\quad
\pi_i^h\circ\pi_j^h=\sum_k b_{ijk}\pi_k^h,\quad b_{ijk}\in\frak m^{h+1},
\]for $i,j=1,\dots, n$ and $i\neq j$.
By commutativity and associativity, one deduces the following constraints on $a_{ij},b_{ijk}$:
\begin{align}
\label{eqid1}
&b_{ijk}=b_{jik},\quad &i,j,k=1,\dots,n,\\
&b_{ijk}\in\frak m^{2h+2},\quad &i,j,k=1,\dots,n,\text{ distinct},\\
\label{eqid3}
&b_{ijj}+a_{ij}\in\frak m^{2h+2},\quad &i,j=1,\dots,n,\quad i\neq j.
\end{align}
Set
\[\pi_i':=\pi_i^h+\sum_j w_{ij}\pi_j^h,\quad i=1,\dots,n,
\]with arbitrary coefficients $w_{ij}\in\frak m^{h+1}$. The $n$-tuple $(\pi_1',\dots, \pi_n')$ is an $(h+1)$-oder idempotent basis\footnote{In fact, the resulting basis $(\pi_1',\dots, \pi_n')$ is not just an $(h+1)$-order idempotent, but even a $(2h+1)$-order idempotent basis.} of $H_K$ if and only if
\[w_{ii}=-a_{ii},\quad \quad
w_{ij}=a_{ij},
\]for $i,j=1,\dots,n$ and $i\neq j$.
This easily follows from \eqref{eqid1}-\eqref{eqid3}.
\endproof
\vskip3mm
In the analytic case, we will say that a \dfm\ $M$ is \emph{(generically) semisimple} if the set $M_{ss}:=\{p\in M\colon (T_pM,\circ_p)\cong \C^n\}$ is non-empty. In such a case, it can be proved that $M_{ss}$ is an open dense subset of $M$. At each point $p\in M_{ss}$ there exists tangent vectors $\pi_1|_p,\dots,\pi_n|_p$ satisfying the relations
\[\pi_i|_p\circ_p\pi_j|_p=\pi_i|_p\delta_{ij},\quad \eta_p(\pi_i|_p,\pi_j|_p)=0,\quad i\neq j.
\]It can be proved that, on sufficiently small open subsets $M_{ss}$, a coherent labeling of the idempotent tangent vectors can be chosen so that the resulting local vector fields are holomorphic. For a detailed discussion and proofs, see \cite[Chapter 2]{Her02}.

\begin{rem}
In both the formal and analytic case we have $e=\sum_i \pi_i$.
\end{rem}

\begin{prop}[\cite{Dub92,Dub96,man}]
For both formal and analytic semisimple Frobenius manifolds, the idempotents vector fields $\pi_1,\dots,\pi_n$ are pairwise commuting, i.e. $[\pi_i,\pi_j]=0$. Equivalently, the dual differential forms $\pi_i^\flat$, defined by $\langle \pi_i^\flat,\pi_j\rangle=\delta_{ij}$, are closed.\qed
\end{prop}

In both the formal and analytic cases, this result implies the existence of a local system of coordinates $\bm u:=(u_1,\dots, u_n)$ such that 
\[du_i=\pi_i^\flat,\quad \frac{\der}{\der u_i}=\pi_i.
\]We will refer to $\bm u$ as the \emph{formal/analytic canonical coordinates}. These functions are defined up to re-ordering and shifts by constants. In the formal case, the functions $u_i$'s are just formal functions, i.e. elements of $k[\![\bm t]\!]$.

\begin{prop}[\cite{Dub92,Dub96,man}]
The formal/analytic canonical coordinates can be uniquely chosen (up to re-ordering) so that $E=\sum_{i=1}^nu_i\frac{\der}{\der u_i}.$\qed
\end{prop}

In all the subsequent part of the paper, we will reserve Latin indices for canonical coordinates $u_1,\dots, n_n$ and their vector fields $\der_i:=\frac{\der}{\der u_i}$. Einstein summation rule will be used only for repeated Greek indices.

\section{Convergence of semisimple formal Frobenius manifolds}\label{sec5}
In this Section we prove the main result of the second part of this paper. 

\begin{thm}\label{thconv}
Let $(H,\Phi,\eta,e,E)$ be a semisimple formal Frobenius manifold over $\C$. Then the domain of convergence of $\Phi$ is non-empty.
\end{thm}

It follows that with any semisimple formal Frobenius manifold there is an associated analytic \dfm, as explained in Section \ref{affa}.
\vskip1mm
For the proof, we require some preliminary material.

\subsection{Extended deformed connection} We introduce one of the main objects attached to Frobenius structures, namely an integrable connection. It can be introduced in both formal and analytic frameworks. 
\vskip3mm
\noindent{\bf Formal case. }Let $k$ be a commutative $\Q$-algebra and $(H,\eta,\Phi)$ a formal Frobenius manifold as in Section \ref{ffm}. Denote by $k(\!(z)\!)$ the $k$-algebra of formal Laurent series in an auxiliary indeterminate $z$. Set $K(\!(z)\!):=k[\![\bm t]\!](\!(z)\!)$ to be the Laurent series with coefficients in $k[\![\bm t]\!]$, and $H_{K(\!(z)\!)}:=H\otimes_k K(\!(z)\!)$. In the following paragraphs we will define two connections on the modules $H_K$ and $H_{K(\!(z)\!)}$ respectively. We firstly recall some basic notions.
\subsubsection{Algebraic connections on modules} Let $A$ be a commutative and unital $k$-algebra, and $P$ an $A$-module. Denote by ${\rm Diff}_1(P,P)$ the set of first order differential operators on $P$, i.e. the $k$-linear morphisms $\mathscr D\in\Hom_k(P,P)$ such that
\[ab\mathscr D(p)-b\mathscr D(ap)-a\mathscr D(bp)+\mathscr D(abp)=0,\quad a,b\in A,\quad p\in P.
\]Both ${\rm Der}_k(A)$ and ${\rm Diff}_1(P,P)$ are naturally equipped with an $A$-module structure. A connection $\nabla$ on $P$ is defined by an $A$-linear morphism $\nabla\colon {\rm Der}_k(A)\to {\rm Diff}_1(P,P)$, $u\mapsto\nabla_u$ satisfying the Leibniz rule
\[\nabla_u(ap)=u(a)p+a\nabla_up,\quad a\in A,\quad p\in P.
\]The curvature of $\nabla$ is the $A$-bilinear morphism $R\colon {\rm Der}_k(A)\times {\rm Der}_k(A)\to \Hom_A(P,P)$ defined by
\[R(u,v):=[\nabla_u,\nabla_v]-\nabla_{[u,v]},\quad u,v\in{\rm Der}_k(A).
\]Given a connection on $P$ we can induce connections on all the tensor products (over $A$) $P^{\otimes p}\otimes \Hom_A(P,A)^{\otimes q}$ by requiring that
\begin{enumerate}
\item $\nabla$ commutes with contractions, 
\item on $A$ (i.e. $p=q=0$) the morphism $\nabla\colon {\rm Der}_k(A)\to {\rm Diff}_1(A,A)$ is just the inclusion.
\end{enumerate}
As a general reference, see e.g. \cite{Sar12}.

\subsubsection{Deformed connections on $H_K$} Consider the case $(A,P)=(K,H_K)$. Define a one-parameter family of connections  $\nabla^{z}\colon {\rm Der}_k(K)\cong H_K\to {\rm Diff}_1(H_K,H_K)$, with $z\in\C$, on the module $H_K$
by the formula
\[\nabla^z_{\Dl_\al}{\Dl_\bt}:=z\Dl_\al\circ \Dl_\bt,\quad \al,\bt=1,\dots, n.
\]
\begin{thm}[\cite{Dub92,man}]\label{ffdc}
WDVV equations \eqref{wdvv} are equivalent to the flatness of $\nabla^z$, for any $z\in \C$.\qed
\end{thm}

\begin{rem}\label{LCrem}
The connection $\nabla:=\nabla^0$ is the  
(formal) Levi-Civita connection for $\eta$, i.e. the unique torsion-free connection satisfying $\nabla\eta=0$. If $(e_1,\dots,e_n)$ is a basis of $H_K$, set $\nabla_{e_i}e_j=\sum_k\Gm_{ij}^ke_k$. One can show that 
\[\Gm_{ij}^k=\frac{1}{2}\sum_\ell \eta^{\ell k}\left(e_i\eta_{jk}+e_j\eta_{ik}-e_k\eta_{ij}\right).
\]The standard differential-geometrical proof works verbatim  in this formal framework. 
\end{rem}
\begin{rem}
The Euler vector field is an affine vector field, i.e. $\nabla\nabla E=0$. This is an equivalent formulation of Lemma \ref{affeu}. 
\end{rem}

\subsubsection{Extended deformed connection on $H_{K(\!(z)\!)}$}\label{fedc} 
We consider now the case $(A,P)=(K(\!(z)\!),$ $H_{K(\!(z)\!)})$.
In what follows we assume that the $K$-linear operator $\nabla^0E\colon {\rm Der}_k(K)\cong H_K\to H_K$ is (diagonalizable and) in diagonal form in the basis $(\Dl_1,\dots,\Dl_n)$. Define two new $K$-linear operators $\mc U,\mu$ by the formulae
\bea
&\mc U\colon H_K\to H_K,\quad &X\mapsto E\circ X,\\
&\mu\colon {\rm Der}_k(K)\cong H_K\to H_K,\quad &X\mapsto \frac{D}{2}-\nabla^0_XE,
\eea
where $D\in k$ is as in \eqref{eu1}. All the tensors $\eta,\circ, \mc U,\mu$ can be $K(\!(z)\!)$-linearly extended to $H_{{K(\!(z)\!)}}$. We will denote such an extension by the same symbols.
\vskip3mm
The extended deformed connection $\widehat\nabla\colon{\rm Der}_k(K(\!(z)\!))\to {\rm Diff}_1(H_{K(\!(z)\!)},H_{K(\!(z)\!)})$ is defined by the formulae 
\[\widehat\nabla_{\frac{\der}{\der t^\al}}X=\nabla^z_{\frac{\der}{\der t^\al}}X,\quad\quad \widehat\nabla_{\frac{\der}{\der z}}X={\frac{\der}{\der z}}X+\mc U(X)-\frac{1}{z}\mu(X),
\]
where $X\in H_{K(\!(z)\!)}$.

\begin{thm}[\cite{Dub96,Dub98,Dub99}]\label{ffedc}
The connection $\widehat\nabla$ is flat.
\end{thm}
\proof The flatness of $\widehat\nabla$ is equivalent to the following conditions:
$\der_\dl\Phi_{\al\bt\gm}$ is completely symmetric in $(\al,\bt,\gm,\dl)$,
the product $\circ$ is associative,
$\nabla\nabla E=0$,
and $\frak L_Ec=c$.
This can be checked by a straightforward computation.
\endproof

\noindent{\bf Analytic case. }Let $M$ be a \dfm. Introduce the $(1,1)$-tensors $\mc U,\mu\in\Gm({\rm End}(TM))$ by the formulae
\[\mc U(X)=E\circ X,\quad \mu(X):=\frac{2-d}{2}X-\nabla_XE,\quad X\in\Gm(TM),
\]where $d$ is the charge of the Dubrovin-Frobenius structure, and $\nabla$ is the Levi-Civita connection of $\eta$. We assume that $\mu$ is (diagonalizable and) in diagonal form in the frame $(\der_{t^1},\dots,\der_{t^n})$. 

Denote by $\pi\colon M\times \C^*\to M$ the canonical projection on the first factor. If $\mathscr T_M$ denotes the tangent sheaf of $M$, then $\pi^*\mathscr T_M$ is the sheaf of sections of $\pi^*TM$, and $\pi^{-1}\mathscr T_M$ is the sheaf of sections of $\pi^*TM$ constant along the fibers of $\pi$.  All the tensors $\eta,c,e,E,\mc U,\mu$ can be lifted to the pulled-back bundle $\pi^*TM$, and we denote these lifts with the same symbols. Consequently, also the Levi-Civita connection $\nabla$ can be uniquely lifted on $\pi^*TM$ in such a way that $\nabla_{\frac{\der}{\der z}}Y=0$ for $Y\in\pi^{-1}\mathscr T_M$.

The extended deformed connection $\widehat\nabla$ is the connection on $\pi^*TM$ defined by the formulae
\beq
\label{aedc1}
\widehat\nabla_{\frac{\der}{\der t^\al}}Y=\nabla_{\frac{\der}{\der t^\al}}Y+z\frac{\der}{\der t^\al}\circ Y,\quad\quad
\widehat\nabla_{\frac{\der}{\der z}}Y=\nabla_{{\frac{\der}{\der z}}}Y+\mc U(Y)-\frac{1}{z}\mu(Y),
\eeq
where $Y\in \pi^*\mathscr T_M$. 
\begin{rem}
If we consider a formal Frobenius manifold associated to a pointed \dfm\ $(M,p)$ as in Section \ref{affa}, the Christoffel symbols of the formal connection $\widehat\nabla$ constructed in Section \ref{fedc} are germs of the Christoffel symbols of \eqref{aedc1}
at the point $p$.
\end{rem}
Theorem \ref{ffedc} and its proof hold verbatim for the connection $\widehat\nabla$ defined by \eqref{aedc1}.

\begin{rem}
In both the formal and analytic case, the operator $\mc U$ is $\eta$-self-adjoint, and $\mu$ is $\eta$-skew-symmetric: for arbitrary $X,Y\in H_K$ (resp. sections of $TM$), we have
\beq
\eta(\mc U(X),Y)=\eta(X,\mc U(Y)),\quad \eta(\mu(X),Y)=- \eta(X,\mu (Y)).
\eeq
\end{rem}

\subsection{Darboux-Egoroff equations} Given a formal (resp.\,\,analytic) semisimple Frobenius manifold with idempotent vectors $\pi_1,\dots,\pi_n$ define the formal (resp.\,\,analytic) functions $\eta_{ii},\gm_{ij}\in k[\![\bm u]\!]$ (resp.\,\,$\C\{\bm u\}$) by
\bea
&\eta_{ii}(\bm u):=\eta(\pi_i(\bm u),\pi_i(\bm u)),\quad &i=1,\dots,n,\\
&\gm_{ij}(\bm u):=\frac{\der_{j}\sqrt{\eta_{ii}(\bm u)}}{\sqrt{\eta_{jj}(\bm u)}},
\quad &i,j=1,\dots,n.
\eea
\begin{lem}We have
\beq\label{symgij}
\gm_{ij}(\bm u)=\frac{1}{2}\frac{\der_i\der_jt_1(\bm u)}{\sqrt{\der_i t_1(\bm u)\der_j t_1(\bm u)}},\quad t_1:=\sum_\al \eta_{1\al}t^\al.
\eeq
In particular, $\gm_{ij}=\gm_{ji}$.
\end{lem}
\proof
Consider the co-unit 1-form $\theta:=\eta(e,-)$.
We have $\eta_{ii}=\eta(\pi_i,\pi_i)=\eta(e,\pi_i)=\langle\theta, \der_i\rangle=\der_it_1$. Equation \eqref{symgij} follows by the definition of $\gm_{ij}$.
\endproof
\begin{thm}\label{tDE}
The functions $\gm_{ij}(\bm u)$ satisfy the Darboux-Egoroff equations, i.e.
\bean
\label{DE1}
&\der_k\gm_{ij}=\gm_{ik}\gm_{kj}, &i,j,k\text{ distinct,}\\
\label{DE2}
&\sum_{k=1}^n\der_k\gm_{ij}=0,\quad &i\neq j\\
\label{DE3}
&\sum_{k=1}^n u_k\der_k\gm_{ij}=-\gm_{ij},\quad &i\neq j.
\eean
\end{thm}
\proof
The proofs of \cite[Prop. 3.4.1, Th. 3.7.2]{man} apply verbatim also to the formal case. Notice that \eqref{DE1} and \eqref{DE2} are equivalent to the flatness of $\eta$.
\endproof

\begin{cor}
For $i\neq j$, we have
\beq
\label{DE4}(u_i-u_j)\der_i\gm_{ij}=\sum_{k\neq i,j}(u_j-u_k)\gm_{ik}\gm_{kj}-\gm_{ij}.
\eeq
\end{cor}
\proof
An easy consequence of \eqref{DE1}, \eqref{DE2}, and \eqref{DE3}.
\endproof

\subsection{$\widehat\nabla$-flatness in canonical coordinates}Given a formal (resp. analytic) semisimple Frobenius manifold with idempotent vectors $\pi_1,\dots,\pi_n$ define the vectors 
\bea
f_i(\bm u):=\eta_{ii}(\bm u)^{-\frac{1}{2}}\pi_i(\bm u),\quad i=1,\dots,n, 
\eea
for some choices of the square roots, and introduce the matrix $\Psi\in GL(n,k[\![\bm u]\!])$ (resp. $GL(n,k\{\bm u\})$) defined by
\[\Psi=(\Psi_{i\al})_{i,\al},\quad \frac{\der}{\der t^\al}=\sum_{i=1}^n \Psi_{i\al}f_i,\quad i=1,\dots,n.
\]

\begin{lem}[{\cite[Lemma 2.29]{CDG20}}]\label{useq}
We have
\[
\pushQED{\qed} 
\Psi^T\Psi=\eta,\quad \Psi_{i1}=\sqrt{\eta_{ii}},\quad \der_i=\sum_{\al,\bt=1}^n\Psi_{i1}\Psi_{i\bt}\eta^{\al\bt}\der_\al,\quad c_{\al\bt\gm}=\sum_{i=1}^n\frac{\Psi_{i\al}\Psi_{i\bt}\Psi_{i\gm}}{\Psi_{i1}}.\qedhere
\popQED
\]
\end{lem}

\begin{lem}\label{manlem}
We have $\mu(f_i)=\sum_{j\neq i}(u_j-u_i)\gm_{ij}f_j$.
\end{lem}
\proof
Set $\nabla_{\pi_i}\pi_j=\sum_k\Gm_{ij}^k\pi_k$. The only nonzero Christoffel symbols are 
\[\Gm_{ii}^i=\frac{1}{2}\eta_{ii}^{-1}\frac{\der\eta_{ii}}{\der u_i},\quad \Gm_{ii}^j=-\frac{1}{2}\eta_{jj}^{-1}\frac{\der\eta_{ii}}{\der u_j},\quad \Gm_{ij}^i=\Gm_{ji}^i=\frac{1}{2}\eta_{ii}^{-1}\frac{\der\eta_{ii}}{\der u_j},\quad i\neq j,
\]
see Remark \ref{LCrem}. The claim follows by straightforward computations.
\endproof

The connection $\widehat\nabla$ can be extended to the full tensor algebra of $H_{K(\!(z)\!)}$ (resp. $\pi^*TM$). 
\begin{lem}Let $\xi\in H^T_{K(\!(z)\!)}$ (resp. $\pi^*T^*M$) such that $\widehat \nabla\xi=0$. If $\xi^\sharp$ is $\eta$-dual to $\xi$, then we have
\begin{align}
\label{sist1fr}
\widehat\nabla_{\frac{\der}{\der t^\al}}\xi^\sharp=2 z\frac{\der}{\der t^\al}\circ \xi^\sharp,
\quad\quad 
\widehat\nabla_{\frac{\der}{\der z}}\xi^\sharp=2\,\mathcal U(\xi^\sharp).
\end{align}
\end{lem}
\proof
A simple computation shows that $\widehat\nabla_{\der_\al}\eta_{\bt\gm}=-2z c_{\al\bt\gm}$ for $\al,\bt,\gm=1,\dots, n$, and $\widehat\nabla_{\der_z}\eta_{\al\bt}=-2\eta_{\la\bt}\mc U^\la_\al$, for $\al,\bt=1,\dots, n$. The claim then follows.
\endproof
Set $\xi^\sharp(\bm t, z)=\sum_{i=1}^n y_i(\bm u(\bm t),z)f_i$ for some formal (resp. analytic) functions $y_i$ of the formal canonical coordinates $\bm u$ and $z$. Then, equations \eqref{sist1fr} are equivalent to
\begin{align}
\label{sist1b}
\frac{\der}{\der u_i}y&=\left(zE_i+V_i(\bm u)\right)y,\\
\label{sist2b}
\frac{\der}{\der z}y&=\left(U(\bm u)+\frac{1}{z}V(\bm u)\right)y, 
\end{align}
where $y=(y_1,\dots, y_n)^T$,
\[V_i:=\frac{\der\Psi}{\der u_i}\Psi^{-1},\quad (E_i)_{ab}=\delta_{ai}\delta_{bi},
\]and
\[U={\rm diag}(u_1,\dots, u_n),\quad V:=\Psi\mu\Psi^{-1},
\]where $\mu$ is the matrix of the operator $\mu\colon H_K\to H_K$ associated with the basis $(\frac{\der}{\der t^1},\dots, \frac{\der}{\der t^n})$.
The compatibility of the system \eqref{sist1b},\eqref{sist2b} is equivalent to the equations
\begin{align}
\label{cmp1}
&\frac{\der V}{\der u_i}=[V_i,V],\\
\label{cmp2}
&[U,V_i]=[E_i,V].
\end{align}

\begin{lem}\label{lvvi}
Set $\Gm=(\gm_{ab})$. We have
\[V^T+V=0,\quad V=[\Gm,U],\quad V^T_i+V_i=0,\quad V_i=[\Gm,E_i],\quad i=1,\dots,n,
\]
\end{lem}
\proof
The identity $V=[\Gm,U]$ is Lemma \ref{manlem}. The identity $\Psi(\bm u)^T\Psi(\bm u)=\eta$ implies $\der_i\Psi^T\Psi+\Psi^T\der_i\Psi=0$, so that $V_i^T+V_i=0$. We have $[U,V_i]=[U,[\Gm,E_i]]$, by \eqref{cmp2} and Jacobi identity. The kernel of the operator $[U,-]\colon M_n(k[\![\bm u]\!])\to M_n(k[\![\bm u]\!])$ consists of diagonal matrices: if $A\in M_n(k[\![\bm u]\!])$ is such that $[U,A]=0$, then $(u_a-u_b)A_{ab}(\bm u)=0$ for any $a,b=1,\dots,n$ with $a\neq b$. We deduce that $A_{ab}(\bm u)=0$, with $a \neq b$, since $k[\![\bm u]\!]$ is an integral domain. Hence $V_i=\mc D+[\Gm,E_i]$, where $\mc D$ is a diagonal matrix. The skew-symmetry of $V_i$ implies that $\mc D=0$.
\endproof

\begin{cor}
The equations
\[V_i=[\Gm,E_i],\quad V=[\Gm, U],\quad \der_i V=[V_i,V],
\]are equivalent to the Darboux-Egoroff equations \eqref{DE1},\eqref{DE2},\eqref{DE3} on $\Gm$.
\end{cor}
\proof
This can be checked by a direct computation.
\endproof

\subsection{Reconstruction of the Frobenius structure}\label{recfrob} 
In this section we recall how it is possible to reconstruct the WDVV potential $\Phi$ of a formal Frobenius manifold starting from the system of $\nabla^z$-flatness for 1-forms.
By Theorem \ref{ffdc} we can look for formal functions $\tilde{\bm t}:=(\tilde t_1,\dots,\tilde t_n)$ of the form
\[\tilde t_\al(\bm t,z):=\sum_{p=0}^\infty h_{\al,p}(\bm t)z^p\in k[\![\bm t,z]\!],\quad h_{\al,0}(\bm t)=t_\al\equiv t^\bt\eta_{\al\bt},
\]such that $\nabla^z d\tilde t_\al=0$ for $\al=1,\dots, n$.

\begin{lem}[\cite{Dub92}]
The functions $h_{\al,p}$ satisfy the recursive equations
\[\pushQED{\qed} 
h_{\al,0}(\bm t)=t_\al\equiv t^\bt\eta_{\al\bt},\quad\quad
\der_{\bt}\der_\gm h_{\al,p+1}=c^\eps_{\bt\gm}\der_\eps h_{\al,p},\quad p\in\N.\qedhere
\popQED
\]

\end{lem}

For $f,g\in K$ write $f\approx g$ iff $f-g$ is a (at most) quadratic polynomial in $\bm t$.
\begin{lem}[\cite{Dub92}]
We have
\bean
\label{hap1} &h_{\al,1}\approx\der_\al\Phi,&\quad \al=1,\dots,n,\\
\label{hap2}
&h_{1,2}\approx t^\al\der_\al\Phi -2\Phi.&
\eean
\end{lem}

\proof
We have $\der_\bt h_{\al,0}=\eta_{\al\bt}$, so that $\der_\gm\der_\bt h_{\al,1}=c_{\al\bt\gm}$. Equation \eqref{hap1} follows.\newline
We have $\der_1\Phi=\frac{1}{2}\eta_{\al\bt}t^\al t^\bt$, so that $\der_\al\der_\bt h_{1,2}=c^\gm_{\al\bt}\der_\gm h_{1,1}=c^\gm_{\al\bt}\der_\gm\der_1\Phi=c^\gm_{\al\bt}\eta_{\gm\nu}t^\nu$. We also have $\der_\al\der_\bt\left(t^\la\der_\la\Phi -2\Phi\right)=c_{\al\bt\la}t^\la$, and \eqref{hap2} follows.
\endproof

Given a function $f\in K$, we denote by $\egr f\in H_K$ the $\eta$-gradient of $f$, defined by $\egr f:=\sum_\al\eta^{\al\bt}\der_\bt f \Dl_\al$. The following result allows to reconstruct the potential $\Phi$ (up to quadratic terms) from the first coefficients $h_{\al,p}$, with $p\leq 3$. 

\begin{thm}[\cite{Dub96,Dub99}]
We have
\beq\label{magfor}
\Phi\approx\frac{1}{2}\left[\eta(\egr h_{\al,1},\egr h_{1,1})\eta^{\al\bt}\eta(\egr h_{\bt,0},\egr h_{1,1})-\eta(\egr h_{1,1},\egr h_{1,2})-\eta(\egr h_{1,3},\egr h_{1,0})\right].
\eeq
\end{thm}

\proof
The expression in square brackets in the r.h.s.\,\,of \eqref{magfor} equals
\bean
\nonumber
&\eta^{\nu\la}\der_\nu h_{\al,1}\der_\la h_{1,1}\eta^{\al\bt}\der_\bt h_{1,1}-\eta^{\tau\eps}\der_\tau h_{1,1}\der_{\eps}h_{1,2}-\der_1 h_{1,3}\\
\label{ap2}
\approx&\eta^{\nu\la}\eta^{\al\bt}\der^2_{\nu\al}\Phi\,\der^2_{\la 1}\Phi\,\der^2_{\bt 1}\Phi-\eta^{\tau\eps}\der^2_{\tau 1}\Phi\left(t^\la\der^2_{\eps\la}\Phi-\der_\eps\Phi\right)-\der_1 h_{1,3}.
\eean
We have $\der^2_{\nu 1}\Phi=\eta_{\nu\al}t^\al$, and $\der_\al\der_1h_{1,3}=c_{1\al}^\bt\der_\bt h_{1,2}=\der_\al h_{1,2}$ so that $\der_1h_{1,3}\approx h_{1,2}$. Hence \eqref{ap2} equals $2\Phi$ up to quadratic terms.
\endproof

\subsection{Proof of Theorem \ref{thconv}}  
 Let $(H,\eta,\Phi, e, E)$ be a formal Frobenius manifold. Fix one ordering $\bm u_o\in\C^n$ of the eigenvalues of $\mc U(\bm t)$ specialized at the origin $\bm t=0$. We have $n\times n$ matrix-valued (a priori) formal power series in $\bm u$
 \bea
 V(\bm u)= V_o+\sum_{k=1}^\infty\sum_{\ell_1,\dots,\ell_k=1}^n \frac{1}{k!}V^{(\bm \ell)}\prod_{j=1}^k\overline{u}_{\ell_j},\quad\quad
 V_i(\bm u)= V_{i,o}+\sum_{k=1}^\infty\sum_{\ell_1,\dots,\ell_k=1}^n \frac{1}{k!}V_i^{(\bm \ell)}\prod_{j=1}^k\overline{u}_{\ell_j},\\
 \Psi(\bm u)= \Psi_o+\sum_{k=1}^\infty\sum_{\ell_1,\dots,\ell_k=1}^n \frac{1}{k!}\Psi^{(\bm \ell)}\prod_{j=1}^k\overline{u}_{\ell_j},\quad\quad
 \Gm(\bm u)= \Gm_o+\sum_{k=1}^\infty\sum_{\ell_1,\dots,\ell_k=1}^n \frac{1}{k!}\Gm^{(\bm \ell)}\prod_{j=1}^k\overline{u}_{\ell_j},
 \eea
 where $\overline{u}_i:=u_i-u_{o,i}$ for $i=1,\dots,n$. These power series are well defined by the semisimplicity assumption, and they satisfy properties described in Theorem \ref{tDE}, and Lemmata \ref{useq}, \ref{manlem} and \ref{lvvi}.
We subdivide the proof in two parts. In the first part, we construct a pointed germ $(M,p)$ of a \dfm\ $M$ starting from the datum of $\bm u_o, V_o, \Psi_o, \Gm_o$. In the second part, we prove that the original formal structure $(H,\eta,\Phi, e, E)$ is the completion of the pointed analytic germ $(M,p)$.
\vskip3mm
\noindent {\bf Part I. }  The system \eqref{sist2b} specialized at $\bm u_o$, namely $\frac{\der Y}{\der z}=(U_o+\frac{1}{z}V_o)Y$, can be identified with equation \eqref{sist1} (in the special case $B_o'=0$). Notice that conditions (1) and (2) of Theorem \ref{tsab} and Proposition \ref{solinf} are satisfied: by Lemma \ref{manlem} we have $V=[\Gm,U]$ and therefore $B'_o=0$ and $B_o''=[\Gm_o, U_o]$. The arguments of Section \ref{ssabt} can be applied, in both cases $\bm u_o\in\C^n\setminus\Dl$ and $\bm u_o\in\Dl$. We can fix an admissible $\tau$ at $\bm u_o$, the $(\bm u_o,\tau)$-admissible datum $\frak M$ is well-defined, and we can set the RHB problem $\mc P[\bm u,\tau,\frak M]$. This problem is solvable with respect to $\bm u$ on an open neighborhood $\mc V\setminus\Theta$ of $\bm u_o$, by Theorem \ref{teoimp}. The unique solution $G(z;\bm u)$ is holomorphic in $\bm u\in\mc V\setminus \Theta$, and with expansion
\bea
&G(z;\bm u)=I+\frac{1}{z}F_1^{\rm an}(\bm u)+O\left(\frac{1}{z^2}\right),&\quad z\to\infty,\quad z\in\Pi_{L/R},\\
&G(z;\bm u)=G_0(\bm u)+G_1(\bm u)z+G_2(\bm u)z^2+G_3(\bm u)z^3+O(z^4),&\quad z\to0.
\eea
Here the superscript ``an'' stands for {\it analytic}.
As output of Section \ref{ssabt}, we also obtain a compatible joint system of differential equations (with analytic coefficients in $\bm u$, not just formal) of the form
\beq\label{ansyst}
\frac{\der Y}{\der u_i}=(zE_i+V_i^{\rm an}(\bm u))Y,\quad \frac{\der Y}{\der z}=\left(U+\frac{1}{z}V^{\rm an}(\bm u)\right)Y,
\eeq
where $V^{\rm an}(\bm u):=[F_1^{\rm an}(\bm u), U],$ and $V_i^{\rm an}(\bm u):=[F_1^{\rm an}(\bm u), E_i]$. Moreover, we have 
\[V^{\rm an}(\bm u_o)=V_o,\quad G_0(\bm u_o)=\Psi_o,\quad \der_i G_0=V_i^{\rm an}G_0,\quad i=1,\dots,n.
\]
From the datum of $G_i(\bm u)$, with $i=0,1,2,3$, we can construct a \dfm\ as follows: set
\bea
t^\al(\bm u)&:=&\eta^{\al\bt}\sum_{i=1}^nG_{0,i\bt}(\bm u)G_{1,i1}(\bm u),\quad \al=1,\dots,n,\\
F(\bm u)&:=&\frac{1}{2}\left[t^\al(\bm u)t^\bt(\bm u)\sum_{i=1}^n G_{0,i\al}(\bm u)G_{1,i\bt}(\bm u)-\sum_{i=1}^n\left(G_{1,i1}(\bm u)G_{2,i1}(\bm u)+G_{0,i1}(\bm u)G_{3,i1}(\bm u)\right)\right].
\eea
Invert the first series expansions, to obtain $\bm u=\bm u(\bm t)$. The function $F(\bm u(\bm t))$ gives a solution of WDVV equations, and defines an analytic \dfm\ on an open subset of $H$. The formulae above are, in their essence, re-writing of formulae of Lemma \ref{useq} and formula \eqref{magfor}. See \cite{Dub99,Guz01}. 
\vskip3mm
\noindent{\bf Part II.} We need to prove that the series expansion $F(\bm u(\bm t))$ obtained in Part I equals (up to quadratic terms) the original potential $\Phi(\bm t)$.
For that, it is sufficient to prove that $F_1^{\rm an}(\bm u)''=\Gm(\bm u)''$. By Lemma \ref{lvvi}, indeed, it follows that $V^{\rm an}(\bm u)=V(\bm u),$ and $V_i^{\rm an}(\bm u)=V_i(\bm u)$. Consequently, from the equations 
\[\der_i G_0=V_i^{\rm an}G_0,\quad \der_i \Psi=V_i\Psi,\quad G_0(\bm u_o)=\Psi_o,
\]we deduce $G_0(\bm u)=\Psi(\bm u)$. From this, the equations defining $F(\bm u(\bm t))$ given in Part I, and the last two formulas of Lemma \ref{useq} (or equivalently the reconstruction formula \eqref{magfor} for $\Phi(\bm t)$), we obtain $\der^3_{\al\bt\gm}\Phi(\bm t)=\der^3_{\al\bt\gm}F(\bm u(\bm t))$. This proves the thesis.
\vskip2mm
We now prove that $F_1^{\rm an}(\bm u)''=\Gm(\bm u)''$.

\begin{lem}
We have $F_1^{\rm an}(\bm u_o)''=\Gm_o''$.
\end{lem}
\proof By Proposition \ref{solinf}, the system \eqref{sist2b} specialized at $\bm u_o$, namely $\frac{\der Y}{\der z}=(U_o+\frac{1}{z}V_o)Y$, admits a unique formal solution $Y_F(z)=\left(I+A_1z^{-1}+A_2z^{-2}+O(z^{-3})\right)e^{zU}$. Let us recall how to compute $A_1$. It is uniquely determined by the two equations
\[[A_1,U_o]=V_o,\quad [A_2,U_o]=A_1+V_oA_1.
\]The first equation uniquely determines all the entries $(A_1)_{ab}$ for indices $a\neq b$ such that $u_{o,a}\neq u_{o,b}$:
\[(A_1)_{ab}=\frac{V_{o,ab}}{u_{o,b}-u_{o,a}}=\Gm_{o,ab},
\]by Lemma \ref{manlem}. All the remaining entries $(A_1)_{ab}$, with $a\neq b$ such that $u_{o,a}=u_{o,b}$, are uniquely determined by the second equation:
\[(A_1)_{ab}=-\sum_{\ell}V_{o,a\ell}(A_1)_{\ell b}=-\sum_{\ell}(u_{o,\ell}-u_{o,a})\Gm_{o,a\ell}\Gm_{\ell b}=\Gm_{o,ab}.
\]The last equality follows by specializing equation \eqref{DE4} to $\bm u=\bm u_o$. This prove that $A_1''=\Gm_o''$. By uniqueness of the formal solution we clearly have $F_1^{\rm an}(\bm u_o)=A_1$.
\endproof

\begin{lem}
The off-diagonal entries of $F_1^{\rm an}(\bm u)$ satisfy the Darboux-Egoroff system \eqref{DE1}, \eqref{DE2}, \eqref{DE3}, \eqref{DE4}.
\end{lem}
\proof From the compatibility conditions $\der_i\der_j=\der_j\der_i$ of the system \eqref{ansyst}, we have
\[[E_j,\der_iF_1^{\rm an}]-[E_i,\der_jF_1^{\rm an}]+[[E_i,F_1^{\rm an}],[E_j,F_1^{\rm an}]]=0.
\]This coincides with equations \eqref{DE1} and \eqref{DE2}. Let $\kappa\in\R_{>0}$. The piecewise analytic function $\widetilde G\colon (\Pi_0\cup\Pi_L\cup\Pi_R)\times(\mc V\setminus\kappa\Theta)\to\C$ defined by
\bea
&\widetilde G(z;\bm u):=G(\kappa z;\kappa^{-1} \bm u)\kappa^Dz^{D}\kappa^L z^{-D},\quad &z\in\Pi_0,\\
& \quad \widetilde G(z;\bm u):= G(\kappa z;\kappa^{-1} \bm u),\quad &z\in\Pi_{L/R},
\eea 
solves the same RHB problem $\mc P[\bm u,\tau,\frak M]$ as $G$. By uniqueness of solution we have $\widetilde G=G$. This implies that $F_1^{\rm an}(\kappa^{-1}\bm u)=\kappa F_1^{\rm an}(\bm u)$, and \eqref{DE3} follows by Euler's homogeneous function theorem.
\endproof

\begin{lem}
Let
\[ \Gm(\bm u)= \Gm_o+\sum_{k=1}^\infty\sum_{\ell_1,\dots,\ell_k=1}^n \frac{1}{k!}\Gm^{(\bm \ell)}\prod_{j=1}^k\overline{u}_{\ell_j},\quad \overline{u}_i:=u_i-u_{o,i},
\]be a matrix-valued formal power series, with $\Gm(\bm u)^T=\Gm(\bm u)$, and whose off-diagonal entries $\Gm_{ij}$ are formal solutions of the Darboux-Egoroff system \eqref{DE1}, \eqref{DE2}, \eqref{DE3}.  
The off-diagonal entries of the coefficients $\Gm^{(\ell)}$ can be uniquely reconstructed from the off-diagonal entries of $\Gm_o$.
\end{lem}
\proof We have to show that the derivatives $\der_{i_1}\dots\der_{i_N}\Gm_{ij}(\bm u_o)$ can be computed from the only knowledge of the numbers $\Gm_{ij}(\bm u_o)$. We proceed by induction on $N$. Let us start with the case $N=1$.\newline
{\bf Step 1.} For $i,j,k$ distinct, by expanding both sides of $\der_k\Gm_{ij}=\Gm_{ik}\Gm_{kj}$ in power series, and equating the coefficients, one reconstructs the coefficients of $\der_k\Gm_{ij}(\bm u_o)$. \newline
{\bf Step 2.} From the identity \eqref{DE4}  for $\Gm_{ij}$, one can compute $\der_i \Gm_{ij}(\bm u_o)$ provided that $u_{o,i}\neq u_{o,j}$.\newline
{\bf Step 3. }Assume that $u_{o,i}=u_{o,j}$. By taking the $\der_i$-derivative of both sides of \eqref{DE4} we obtain
\beq
\label{diDE4}
2\der_{i}\Gm_{ij}(\bm u)+(u_i-u_j)\der_i\der_i \Gm_{ij}(\bm u)=\sum_{k\neq i,j}(u_j-u_k)\left[\der_i \Gm_{ik}(\bm u) \Gm_{kj}(\bm u)+ \Gm_{ik}(\bm u)\der_i \Gm_{kj}(\bm u)\right].
\eeq
By evaluating \eqref{diDE4} at $\bm u=\bm u_o$ we can compute all the numbers $\der_i \Gm_{ij}(\bm u_o)$, namely 
\[\der_i \Gm_{ij}(\bm u_o)=\frac{1}{2}\sum_{k\neq i,j}(u_{o,j}-u_{o,k})\left[\der_i \Gm_{ik}(\bm u_o) \Gm_{kj}(\bm u_o)+ \Gm_{ik}(\bm u_o)^2 \Gm_{ij}(\bm u_o)\right].
\]Notice that the only terms $\der_i \Gm_{ik}(\bm u_o)$ appearing in this sum are those computed in Step 2.\newline
{\bf Step 4.} By the symmetry condition $\Gm(\bm u)^T=\Gm(\bm u)$, we have $\der_j\Gm_{ij}(\bm u_o)=\der_j\Gm_{ji}(\bm u_o)$, and these numbers can be computed as in Steps 2 and 3.\newline
This proves that all the first derivatives $\der_k \Gm_{ij}(\bm u_o)$ can be computed.
\vskip2mm
\noindent{\bf Inductive step. }Assume to know all the $N$-th derivatives $\der_{i_1}\dots \der_{i_N} \Gm_{ij}(\bm u_o)$. We show how to compute the number $\der_{h_1}\dots \der_{h_{N+1}} \Gm_{ij}(\bm u_o)$ for any $(N+1)$-tuple $(h_1,\dots, h_{N+1})$.\newline 
{\bf Step 1.} Assume that there exists $\ell\in\{1,\dots, N+1\}$ such that $h_\ell\neq i,j$. We have
\[\der_{h_1}\dots\der_{h_{N+1}}\Gm_{ij}=\der_{h_1}\dots\der_{h_{\ell-1}}\der_{h_{\ell+1}}\dots\der_{h_{N+1}}[\der_{h_\ell}\Gm_{ij}]=\der_{h_1}\dots\der_{h_{\ell-1}}\der_{h_{\ell+1}}\dots\der_{h_{N+1}}[\Gm_{ih_{\ell}}\Gm_{h_{\ell}j}].
\]By evaluation at $\bm u=\bm u_o$, we can compute all the numbers $\der_{h_1}\dots\der_{h_{N+1}}\Gm_{ij}(\bm u_o)$.
\vskip2mm
\noindent Now we need to compute the mixed derivatives $\der_i^p\der_j^{N+1-p}\Gm_{ij}(\bm u_o)$, with $0\leq p\leq N+1$. \newline
{\bf Step 2.} Assume $p>0$ and $u_{o,i}\neq u_{o,j}$. Take the $\der_i^{p-1}\der_j^{N+1-p}$-derivative of both sides of \eqref{DE4}: by evaluation at $\bm u=\bm u_o$ we can reconstruct the numbers $\der_i^{p}\der_j^{N+1-p}\Gm_{ij}(\bm u_o)$.\newline
{\bf Step 3.} Assume $p>0$ and $u_{o,i}= u_{o,j}$. Take the $\der_i^{p-1}\der_j^{N+1-p}$-derivative of both sides of \eqref{diDE4}, to obtain
\begin{multline}\label{lisboa1}(p+1)\der_i^p\der_j^{N+1-p}\Gm_{ij}-(N+1-p)\der_i^{p+1}\der_j^{N-p}\Gm_{ij}+(u_i-u_j)\der_i^{p+1}\der_j^{N+1-p}\Gm_{ij}\\
=\der_i^{p-1}\der_j^{N+1-p}\sum_{k\neq i,j}(u_j-u_k)[\der_i \Gm_{ik} \Gm_{kj}+\Gm_{ik}^2\Gm_{ij}].
\end{multline}
Specialize \eqref{lisboa1} for $p=N+1$: by evaluation at $\bm u=\bm u_o$ of both sides, we can compute the derivative $\der_i^{N+1}\Gm_{ij}(\bm u_o)$.\newline
Specialize \eqref{lisboa1} for $p=N$: by evaluation at $\bm u=\bm u_o$ of both sides, we can compute the derivative $\der_i^{N}\der_j\Gm_{ij}(\bm u_o)$.\newline
Repeating this procedure, by decreasing $p\mapsto p-1$ at each step, we can compute all the mixed derivatives $\der_i^p\der_j^{N+1-p}\Gm_{ij}(\bm u_o)$.\newline
{\bf Step 4.} Assume $p=0$. By symmetry of $\Gm(\bm u)$, we have $\der_j^{N+1}\Gm_{ij}(\bm u_o)=\der_j^{N+1}\Gm_{ji}(\bm u_o)$, and we can proceed as in Steps 2 and 3.
\vskip2mm
\noindent This proves that all the $(N+1)$-th derivatives $\der_{h_1}\dots\der_{h_{N+1}}\Gm_{ij}(\bm u_0)$ can be computed.
\endproof

This proves that $F_1^{\rm an}(\bm u)''=\Gm(\bm u)''$. The proof of Theorem \ref{thconv} is complete.

\section{Application to CohFT's and Gromov-Witten theory}\label{sec6}
\subsection{Cohomological field theories} Let $k$ and $(H,\eta,e)$ be as in Section \ref{ffm}. For a pair of nonnegative integers $(g, \frak n)$ in the stable range $2g-2+\frak n > 0$, denote by $\overline{\mc M}_{g,\frak n}$ the Deligne-Mumford moduli space of stable $\frak n$-pointed curves of genus $g$. Denote by $\pi\colon \overline{\mc M}_{g,\frak n+1}\to \overline{\mc M}_{g,\frak n}$ the morphism forgetting the last puncture, by $\si\colon \overline{\mc M}_{g_1,\frak n_1+1}\times \overline{\mc M}_{g_2,\frak n_2+1}\to \overline{\mc M}_{g_1+g_2,\frak n_1+\frak n_2}$ the morphism which identifies the last markings, and by $\tau\colon \overline{\mc M}_{g,\frak n+2}\to \overline{\mc M}_{g+1,\frak n}$ the morphism identifying the last two punctures of a same curve.
\vskip3mm
A {\it Cohomological field theory} (CohFT) on $(H,\eta,e)$ is the datum of a system $(\Om_{g,\frak n})_{2g-2+\frak n>0}$ of $k$-multilinear maps $\Om_{g,\frak n}\colon H^{\otimes \frak n}\to H^\bullet(\overline{\mc M}_{g,\frak n},k)$ satisfying the following axioms:
\begin{enumerate}
\item each tensor $\Om_{g,\frak n}$ is $\frak S_\frak n$-covariant with respect to the natural actions of the symmetric group $\frak S_\frak n$ on both $H^{\otimes \frak n}$ and $H^\bullet(\overline{\mc M}_{g,\frak n},k)$,
\item $\Om_{0,3}(e\otimes\Dl_\al\otimes\Dl_\bt)=\eta_{\al\bt}$,
\item $\pi^*\Om_{g,\frak n}({\bigotimes}_{i=1}^\frak n v_{\al_i})=\Om_{g,\frak n}({\bigotimes}_{i=1}^\frak n v_{\al_i}\otimes e)$,
\item $\sigma^*\Om_{g_1+g_2,\frak n_1+\frak n_2}(\bigotimes_{i=1}^{\frak n_1+\frak n_2}v_{\al_i})=\eta^{\mu\nu}\Om_{g_1,\frak n_1+1}(\bigotimes_{i=1}^{\frak n_1}v_{\al_i}\otimes \Dl_\mu)\Om_{g_2,\frak n_2+1}(\bigotimes_{i=\frak n_1+1}^{\frak n_1+\frak n_2}v_{\al_i}\otimes \Dl_\nu)$,
\item $\tau^*\Om_{g+1,\frak n}(\bigotimes_{i=1}^{\frak n}v_{\al_i})=\eta^{\mu\nu}\Om_{g,\frak n+2}(\bigotimes_{i=1}^{\frak n}v_{\al_i}\otimes\Dl_\mu\otimes\Dl_\nu)$.
\end{enumerate}
Given a CohFT, we may introduce generating functions, in infinitely many variables $\bm t^\bullet_\bullet=(t^\al_d)_{\substack{\al=1,\dots, n\\ d\in\N}}$ of intersection numbers with psi-classes,
\bean
\label{despot}
\mc F_{g}( \bm t^\bullet_\bullet):=\sum_{\substack{\frak n\geq 0\\ 2g-2+\frak n>0}}\frac{1}{\frak n!}\sum_{\substack{\al_1,\dots,\al_\frak n=1,\dots,n \\ d_1,\dots, d_\frak n\geq 0}}\Big\langle\prod_{i=1}^{\frak n}\tau_{d_i}\Dl_{\al_i}\Big\rangle_g\prod_{i=1}^{\frak n}t^{\al_i}_{d_i},\\
\Big\langle\prod_{i=1}^{\frak n}\tau_{d_i}\Dl_{\al_i}\Big\rangle_g:=\int_{\overline{\mc M}_{g,\frak n}}\Om_{g,\frak n}\left(\bigotimes_{i=1}^\frak n \Dl_{\al_i}\right)\prod_{i=1}^\frak n\psi_i^{d_i}.
\eean
In the genus zero sector and restricting to the small phase space, i.e. by setting $t^\al_d=0$ for $d>0$ and $t^\al_{0}=t^{\al}$ for $\al=1,\dots,n$, the expression above simplifies to 
\beq
\label{potcohft}
\mc F_0(\bm t)=\sum_{\substack{\frak n>2}}\sum_{\substack{\al_1,\dots,\al_\frak n=1}}^n\frac{t^{\al_1}\dots t^{\al_\frak n}}{\frak n!}\int_{\overline{\mc M}_{0,\frak n}}\Om_{0,\frak n}\left(\Dl_{\al_1}\otimes\dots\otimes\Dl_{\al_\frak n}\right).
\eeq
The power series $\mc F_0\in k[\![\bm t]\!]$ is a solution of WDVV equations, and it defines a formal Frobenius manifold (over $k$) on $(H,\eta,e)$, see \cite{KM94,man}. The CohFT will be said to be \emph{semisimple} if the corresponding formal Frobenius manifold is semisimple. 

If $E=\sum_\al\left(w_\al t^\al+y_\al\right)\der_\al$  is a Killing-conformal vector field on $H$, i.e. $\frak L_E\eta=(2-d)\eta$ for some $d\in k$, we have a natural action of $E$ on the CohFT $(\Om_{g,\frak n})_{g,\frak n}$. Denote by ${\rm deg}\colon H^\bullet(\overline{\mc M}_{g,\frak n},k)\to H^\bullet(\overline{\mc M}_{g,\frak n},k)$ the operator which acts on $H^{2k}$ by multiplication by $k$. Then we set
\bea
(E\Om)_{g,\frak n}\left(\bigotimes_{j=1}^{\frak n} \Dl_{\al_j}\right):=\left({\rm deg}+\sum_{\ell=1}^{\frak n}w_\ell\right)\Om_{g,\frak n}\left(\bigotimes_{j=1}^{\frak n} \Dl_{\al_j}\right)+\pi_*\Om_{g,\frak n+1}\left(\bigotimes_{j=1}^{\frak n} \Dl_{\al_j}\otimes\sum_{\ell=1}^{\frak n}y_\ell\Dl_\ell\right).
\eea
A CohFT is called \emph{ homogeneous in genus $g$} if $(E\Om)_{g,\frak n}=[(g-1)d+n]\Om_{g,\frak n}$ for all $\frak n>2-2g$. When a CohFT is homogeneous in genus zero, $E$ is an Euler vector field for the underlying formal Frobenius manifold. 

\begin{rem}
Teleman Reconstruction Theorem \cite[Th.\,\,1]{Tel12} asserts that a CohFT, semisimple and homogeneous in all genera, can be uniquely reconstructed from the underlying formal Frobenius manifold. The reconstruction is performed via the Givental group action \cite{Giv01}. 
\end{rem}
The following result immediately follows from Theorem \ref{thconv}.
\begin{thm}
For any semisimple and homogeneous (at least in genus 0) CohFT over $k=\C$, the potential $\mc F_0(\bm t)$ is convergent.
In particular, there exist real positive constants $m,\rho_1,\dots,\rho_n$ such that
\[\left|\int_{\overline{\mc M}_{0,|\bm \al|}}\Om_{0,|\bm \al|}\left(\Dl_1^{\otimes \al_1}\otimes\dots\otimes\Dl_n^{\otimes \al_n}\right)\right|\leq m\,\bm\al!\prod_{i=1}^n\rho_i^{\al_i},\quad \bm\al\in\N^n,
\]where we set $\bm\al!:=\prod_{j}\al_j$, and $|\bm \al|:=\sum_{k}\al_k$.\qed
\end{thm}

\subsection{Gromov-Witten theory}Let $X$ be a smooth complex projective variety with vanishing odd cohomology $H^{\rm odd}(X;\C)=0$. Let $(\Dl_1,\dots,\Dl_n)$ be a homogeneous basis of $H^\bullet(X;\C)$, with $\Dl_1=1$ and $(\Dl_2,\dots,\Dl_{r+1})$ a NEF\footnote{This means that each $\Dl_2,\dots,\Dl_{r+1}$ intersect every effective curve class $\bt\in{\rm Eff(X)}$ non-negatively.} $\Z$-basis of $H^2(X;\Z)/$torsion. Denote by $\eta$ the Poincar\'e metric $\eta(\al,\bt):=\int_X\al\cup \bt$. Introduce indeterminates ${\bf Q}:=(Q_1,\dots, Q_r)$, and define the {\it Novikov ring} $\La:=\Q[\![\Qf]\!]$. 
\vskip2mm
Gromov-Witten theory naturally provides a CohFT over the $\La$-module $H^\bullet(X;\La)$ with $\La$-bilinearly extended Poincar\'e metric $\eta$. The maps $\Om_{g,\frak n}$ are given by the counting of curves on $X$,
\beq\label{gwom}
\Om_{g,\frak n}\left(\bigotimes_{i=1}^\frak n\Dl_{\al_i}\right):=\sum_\bt \pho_*\left([\overline{\mc M}_{g,\frak n}(X,\bt)]^{{\rm vir}}\cap\prod_{i=1}^\frak n{\rm ev}_i^*\Dl_{\al_i}\right)\Qf^\bt\in H^\bullet(\overline{\mc M}_{0,\frak n};\La),
\eeq where $\Qf^\bt:=\prod_{i=1}^r Q_i^{\int_\bt \Dl_{i+1}}$, $\overline{\mc M}_{g,\frak n}(X,\bt)$ is the Deligne-Mumford moduli space of $\frak n$-pointed stable maps with target $X$, genus $g$ and degree $\bt$, ${\rm ev}_i\colon \overline{\mc M}_{g,\frak n}(X,\bt)\to X$ are the evaluation morphisms and $\pho\colon \overline{\mc M}_{g,\frak n}(X,\bt)\to \overline{\mc M}_{g,\frak n}$ is the morphism forgetting the map.
\vskip2mm
Equation \eqref{potcohft} defines then a formal power series $F^X_0\in\La[\![\bm t]\!]$, called the genus 0 Gromov-Witten potential of $X$. The corresponding formal Frobenius manifold over $k=\La$ is the \emph{quantum cohomology} of $X$. In order to work with formal Frobenius manifolds over $\C$ we make the following assumption.
\vskip2mm
\noindent{\bf Assumption A:} There exist a point $\bm q\in\C^r$ such that the series $\int_{\overline{\mc M}_{0,\frak n}}\left.\Om_{0,\frak n}\left(\bigotimes_{i=1}^\frak n\Dl_{\al_i}\right)\right|_{\Qf=\bm q}$ are convergent for any $\frak n\geq 3$.
\vskip2mm
If Assumption A holds true, then the specialization $F^X_0|_{\Qf=\bm q}$ is a formal power series in $\C[\![\bm t]\!]$. We call \emph{big quantum cohomology of $X$ }(\emph{at $\Qf=\bm q$}) the corresponding 
formal Frobenius ma\-nifold over $\C$.
We call \emph{small quantum cohomology of $X$ }(\emph{at $\Qf=\bm q$}) the Frobenius $\C$-algebra structure defined on $H^\bullet(X;\C)$ with structure constants $c^\gm_{\al\bt}:=\eta^{\gm\mu}\int_{\overline{\mc M}_{0, 3}}\left.\Om_{0, 3}\left(\Dl_\al\Dl_\bt\Dl_\mu\right)\right|_{\Qf=\bm q}$.

\begin{rem}
Assumption A holds true for all Fano varieties. 
This is because any sum $\sum_\bt$ in \eqref{gwom} reduces to a finite number of terms, so that $\Om_{0,n}\left(\bigotimes_{i=1}^\frak n\Dl_{\al_i}\right)\in\Q[\Qf]$. See e.g. \cite[Prop. 8.1.3]{CK99}.
\end{rem}

\begin{rem}
By the Divisor axiom of Gromov-Witten invariants, it follows that the potential $F^X_0$ can be seen as a formal power series in $\Q[\![t^1,Q_1e^{t^2},\dots, Q_r e^{t^{r+1}}, t^{r+2},\dots, t^n]\!]$, see \cite{CK99,man}. If Assumption A holds true, without loss of generalities we can assume that $\bm q=(1,1,1,\dots,1)$: this correspond to a shift of coordinates $t^{i+1}\mapsto t^{i+1}-\log q_i$ for $i=1,\dots, r$.
\end{rem}

\begin{rem}
If $X$ has generically semisimple quantum cohomology (as a formal Frobenius manifold over $\La$), then $X$ is of Hodge-Tate type, i.e. the Hodge numbers $h_{p,q}(X):=\dim_\C H^q(X,\Om^p)$ vanish for $p\neq q$, see \cite{HMT}.
\end{rem}

Theorem \ref{thconv} implies then the following result.
\begin{thm}\label{thgw}
Let Assumption A hold true. If the small quantum cohomology of $X$ at $\bm q$ is semisimple, then the function $F^X_0(\bm t)|_{\Qf=\bm q}$ has a non-empty domain of convergence $M_{\bm q}\subseteq H^\bullet(X;\C)$, which is equipped with a \dfm\ structure.\qed
\end{thm}

Theorem \ref{thgw} should be compared with other results in literature, differing in techniques. 
In \cite{Iri07}, H.\,Iritani proved convergence of the big quantum cohomology of $X$ under a different assumption, namely that $H^\bullet(X;\C)$ is generated by $H^2(X;\C)$, 
see \cite[Corollary 5.9]{Iri07}. Subsequently, in \cite{CI15} T.\,Coates and H.\,Iritani proved the convergence (suitably defined) of all potentials $\mc F_g^X$ given by \eqref{despot}, by assuming both convergence of $F^X_0$ and semisimplicity.

Whenever the three-point Gromov-Witten correlators $\int_{\overline{\mc M}_{0, 3}}\Om_{0, 3}\left(\Dl_\al\Dl_\bt\Dl_\mu\right)$ of $X$ are explicitly 
known, and thus generators and relations for the small quantum cohomology
ring are given, it is purely a problem in computational commutative algebra to check generic semisimplicity of the small quantum cohomology. Here, we limit ourselves to the following claim\footnote{Surely enough, such a list does not cover all the known cases of semisimple small quantum cohomologies available in literature.}, which  follows from \cite{BM04,BM19,CMP10,Cio,Cio05,Iri07,Per}.

\begin{cor}
We have $F^X_0\in\Q\{\Qf,\bm t\}$ in the following cases (not mutually excluding):
\begin{enumerate}
\item $X=G/P$ is a (co)minuscule homogeneous variety;
\item $X$ is a del Pezzo surface;
\item $X$ is a Fano toric variety;
\item $X$ is one of the following Fano threefolds:
\begin{itemize}
\item $\Pb^3$, a quadric $Q_3,V_5,V_{22}$, 
\item $M^2_k$ with $21\leq k\leq 36$ and $k\neq 23,25,28$, 
\item $M^3_k$ with $k=10, 12, 15, 17, 18, 20, 24, 25, 27, 28, 30, 31$,
\item $\Pb^1\times \Pb^2_k$ where $\Pb^2_k$ is the blow-up of $\Pb^2$ at $k$ points ($1 \leq k\leq 8$);
\end{itemize}
\item $X$ is a Fano general hyperplane section with index $i(X)>\frac{1}{2}{\dim_\C X}$ of a homogeneous space in the following list:
\bea
\Pb^n,\quad \text{the $n$-dimensional quadric $Q_n$,}\quad LG(3,6),\quad F_4/P_1\\
Gr(2,2n+1),\quad OG(5,10),\quad OG(2,2n+1),\quad G_2/P_1;
\eea
\item $X$ is the Cayley Grassmannian parametrizing four dimensional
subalgebras of the complex octonions. \qed
\end{enumerate}
\end{cor}
\begin{rem}
It is known that there exist homogeneous spaces with non-semisimple small quantum cohomology, \cite{CMP10,CP11}. Isotropic Grassmannians $IG(2,2n)$ furnish an example. It is also known, however, that their big quantum cohomology is generically semisimple \cite{GMS15,Per,CMMPS19}. For these varieties, the results of the current paper do not allow to infer the convergence of the genus zero Gromov--Witten potential, a working assumption in \cite[Th. B]{CMMPS19}.
\end{rem}

There is an intriguing conjecture due to B.\,Dubrovin \cite[Conj.\,\,4.2.2]{Dub98} stating the equivalence of the semisimplicity of the (big) quantum cohomology of a variety $X$ (originally assumed to be Fano) and the existence of full exceptional collections in the derived category of coherent sheaves $\mc D^b(X)$. In its most updated formulation, under the assumption of convergence of the genus zero Gromov-Witten potential $F^X_0$, Dubrovin's conjecture also predicts the monodromy data of the system \eqref{sist2b} (in the terminology of the current paper, the admissible data $\frak M$) in terms of characteristic classes of the objects of these exceptional collections, see \cite{GGI16,CDG18,Cot20}. In \cite[\S 4.2, Problem 1]{Dub98} Dubrovin also briefly addressed the problem of convergence of the genus zero Gromov-Witten potential $F^X_0$. In this regard, Dubrovin adds: \guillemotleft {\it Hopefully, in the semisimple case the convergence can be proved on the basis
of the differential equations of n.3\,}\guillemotright\footnote{These include the equations in Part I of our proof of Theorem \ref{thconv}.}. Theorems \ref{thgw} fulfills Dubrovin's hope.
\vskip2mm

\end{document}